\documentclass[10pt]{scrartcl}

\usepackage{amsmath,amsthm,amsfonts,amssymb,amscd}
\usepackage{subfigure}
\usepackage{graphicx,bm,color}
\usepackage{algorithm}
\usepackage{algpseudocode}
\usepackage{stmaryrd}
\usepackage{booktabs,ctable,multirow,longtable}

\usepackage{hyperref}

\definecolor{gray}{gray}{0.6}

\newtheorem{Remark}{Remark}[section]

\newtheorem{form}{Formulation}[section]

\def\ep{\varepsilon}
\def\eps{\varepsilon}
\def\R{\mathbb R}

\begin{document}

\title{A Phase-Field Description for 
 Pressurized and Non-Isothermal Propagating Fractures}

\author{Nima Noii and Thomas Wick\\[1em]
Leibniz Universit\"at Hannover\\
Institut f\"ur Angewandte Mathematik\\
AG Wissenschaftliches Rechnen\\
Welfengarten 1, 30167 Hannover, Germany
\thanks{email: noii@ifam.uni-hannover.de, thomas.wick@ifam.uni-hannover.de}
}


\date{{\small
Dec 22, 2018\\
This is the preprint version (first submission) of an accepted article 
to be published in\\ Computer Methods in Applied Mechanics and Engineering
(CMAME)
\url{https://www.journals.elsevier.com/computer-methods-in-applied-mechanics-and-engineering}
%
}
}

\maketitle

\begin{abstract}
In this work, we extend a phase-field approach for pressurized fractures to non-isothermal settings.
Specifically, the pressure and the temperature are given quantities and the emphasis is on the correct 
modeling of the interface laws between a porous medium and the fracture. The
resulting model is augmented with thermodynamical arguments and then analyzed
from a mechanical perspective. The numerical solution is based on a robust
semi-smooth Newton approach in which the linear equation systems are solved
with a generalized minimal residual method and algebraic multigrid
preconditioning. The proposed modeling and algorithmic developments are
substantiated with different examples in two- and three dimensions. We notice
that for some of these configurations manufactured solutions can be
constructed, allowing for a careful verification of our
implementation. Furthermore, crack-oriented predictor-corrector adaptivity and
a parallel implementation are used to keep the computational cost
reasonable. Snapshots of iteration numbers show an excellent performance of
the nonlinear and linear solution algorithms.
Lastly, for some tests, a computational analysis of the effects of 
strain-energy splitting is performed, which has not been undertaken to date
for similar phase-field settings involving pressure, fluids or 
non-isothermal effects.
\end{abstract}

\begin{keyword}
phase-field; pressurized and non-isothermal fractures; interface laws;
finite elements; adaptivity; parallel programming code
\end{keyword}

%
%
%

\section{Introduction} 
\label{sec:Intro}
In subsurface hydraulic fractures, it is known that cold injections (e.g.,
fluids or $CO2$) may open and advance 
fractures \cite{Ha81,WaPa99,GhaKu06,Chun13,TrSeNg13,GhoAzRa15,SALIMZADEH2018130,SALIMZADEH2018212,doi:10.1029/GM042p0019}.
These effects play a role in oil/gas production,  $CO2$ storage, and
geothermal energy production. For instance, 
\cite{GhaKu06} investigated the effects of cold-water injection on the 
growth and decay of a pre-existing crack when pressure and 
temperature are taken into account. They developed an analytical
formula to solve their model.
Another study incorporating thermal effects was conducted in
\cite{Chun13}. Therein, the fracture propagation is 
based on the displacement discontinuity method. 
For similar computations the displacement 
discontinuity method was used in \cite{GhaZha2006}.
In \cite{WaPa99}, the interaction of temperature and fluid stresses 
with respect to wellbores have been investigated. The authors show
that warm fluids avoid (or restrict) fracture development.
Recently, \cite{SALIMZADEH2018130,SALIMZADEH2018212} studied 
$CO2$ injections and their influence on the fracture aperture/growth 
and a thermo-hydro-mechanical 
approach for fractured geothermal systems.

In the previously mentioned studies, different 
numerical discretization techniques were adopted.
A popular variational approach to fracture 
was introduced in \cite{FraMar98} (with 
according numerics presented in \cite{BourFraMar00}) 
and which was complemented with thermodynamical arguments 
in \cite{AmorMarigoMaurini2009} and \cite{MieWelHof10,MieWelHof10b}.
In the last study, the alternative name phase-field fracture 
was given. The extension of phase-field to hydraulic fractures in poroelasticity 
was first undertaken in \cite{MiWheWi18} (originally published as ICES preprint \cite{MiWheWi13a}).
Therein the solution of the fluid equation (of Darcy type) was assumed to be given, yielding 
a given pressure for the poroelastic mechanics step. On the other hand,
phase-field models involving thermal effects in pure elasticity or plasticity 
have been studied in \cite{BouMaMauSi14,Miehe2015449,Miehe2015486}.

The setting presented in \cite{MiWheWi18} constitutes the starting point to formulate 
the mechanics step of a \textit{phase-field fracture model in
  thermo-poroelasticity} with pressure and temperature as given quantities.
For the temperature contribution, an elliptic diffusion equation
is considered from which a time-dependent decline 
constant can be derived; see Hagoort \cite{Ha81}.
This parameter controls 
the thermo-poroelastic stress (also called back-stress). 
Based on analytical formulas
derived by Sneddon and Lowengrub \cite{Snedd46,SneddLow69}
for the crack opening displacements (COD,
aperture) for pressurized fractures in an elastic medium,
Hagoort extended these equations to non-isothermal 
situations \cite{Ha81}. A very detailed derivation and discussions 
of the application in thermo-poroelasticity can be found in Tran et al. \cite{TrSeNg13}.

In this paper, we follow the ideas outlined in \cite{TrSeNg13}
and develop a phase-field model for 
pressurized and non-isothermal fractures. 
Here, we use the ideas from \cite{MiWheWi18} and carefully 
derive interface conditions between the thermo-poroelastic medium 
and the fracture. The resulting models (i.e., the energy functional and the 
related Euler-Lagrange equations) are then analyzed in detail 
from a mechanical perspective. Here, the focus is on 
thermodynamical arguments and strain-energy splitting, 
which are also complemented with corresponding numerical tests.

Our proposed model is implemented into the adaptive 
parallel framework developed in \cite{HeWheWi15,WiLeeWhe15_eccomas,LeeWheWi16}
with most recent results on its computational performance provided in \cite{HeiWi18}.
A key purpose in the current work is on a detailed code verification with 
respect to the chosen discretization, robustness and efficiency.  
Therefore, some of our settings (in 2D and 3D) are compared with 
manufactured solutions developed in \cite{TrSeNg13}. Moreover,
in some tests, a computational analysis of the effects of 
strain-energy splitting is performed in order to study the influence 
in pressurized and non-isothermal fracture propagation. 
Iteration numbers of the linear solver and the evolution of 
the nonlinear residual and the primal-dual active set (enforcing 
the crack irreversibility constraint) method are undertaken to 
study the performance of the numerical algorithms.

The outline of this paper is as follows: In Section \ref{sec_noniso}
a phase-field model accounting for given pressures and given temperatures 
is derived. Then, in Section \ref{sec_thermo} the model is augmented 
with thermodynamical arguments and then analyzed from a mechanical
perspective. The final model and the numerical solution are discussed 
in Section \ref{sec_numerics}. Several tests demonstrating 
our developed model are presented in Section \ref{sec_tests}.
Here, comparisons to analytical solutions, grid refinement studies, 
calculation of quantities of interest such as the aperture, adaptive 
studies and the performance of the nonlinear and linear solvers are provided.
In Section \ref{sec_conc} our main findings are summarized.

\newpage
\section{Non-isothermal and pressurized phase-field modeling}
\label{sec_noniso}
In the following,
let $\Omega \in \R^d$, $d=2,3$ be a smooth open and bounded set. In $\Omega$, a lower dimensional fracture is denoted by $\mathcal{C}\in \R^{d-1}$. We assume either Dirichlet boundaries conditions $\partial\Omega_D :=
\partial\Omega$ and Neumann condition on $\partial_N \Omega := \Gamma_N \cup \partial\mathcal{C}$, where $\Gamma_N$ 
denotes the outer domain boundary and $\partial\mathcal{C}$ the crack boundary. Let $I:=(0,T)$ denote the loading/time interval with $T>0$ being the end time value.

\begin{figure}[!ht]
	\centering
	{\includegraphics[clip,trim=0cm 3cm 0cm 5cm, width=15cm]{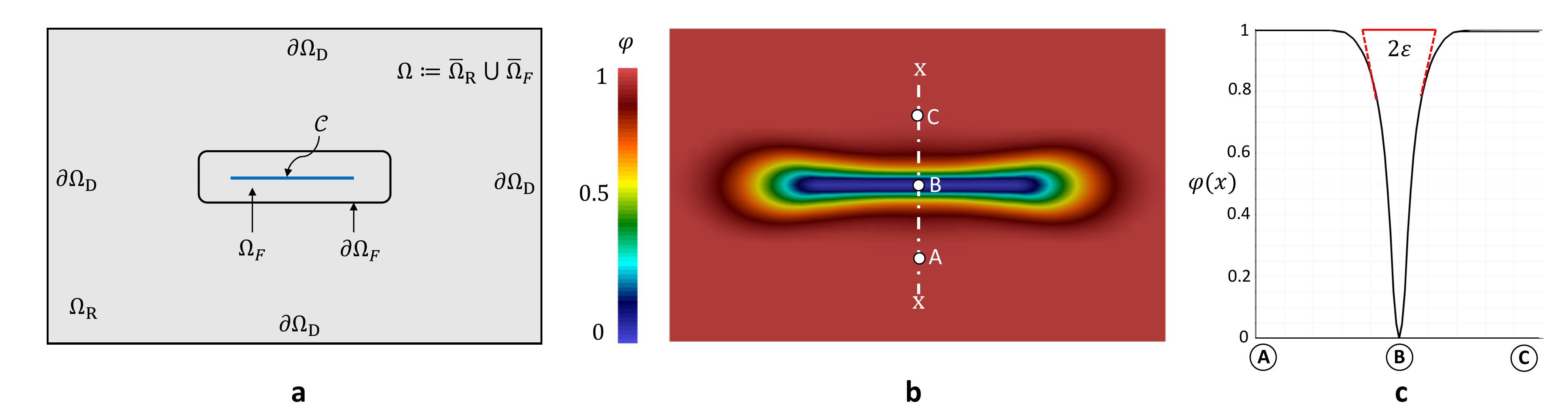}}  
	\caption{(a) Setup of the notation: the unbroken domain is denoted by
          $\Omega$ and $\mathcal{C}$ is the crack phase-field. The smeared crack
          phase-field is approximated by the domain $\Omega_F$. The half
          thickness of $\Omega_F$ is $\ep$. The fracture boundary is $\partial
          \Omega_F$ and the outer boundary of the domain is $\partial
          \Omega$.(b) The corresponding realization using phase-field is shown. Here, the lower-dimensional fracture
          ($\varphi=0$) is approximated with the phase-field
          variable. Consequently, $\Omega_F$ can be represented in terms of
          $\varphi$. (c) The transition zone with $0< \varphi< 1$ has the
          thickness of $2\ep$ where the crack profile on the section $x-x$
          shown in (b) is obtained through a finite element method for solving
          the multi-field (i.e. $\bm u \text{ and } \varphi$) problem.
	}
	\label{Figure1}
\end{figure}

Prototype configurations of the setting are illustrated in Fig. \ref{Figure1}. Using a phase-field approach, the surface fracture $\mathcal{C}$ is approximated in $\Omega_F\subset\Omega \in \mathbb{R}^d$. The intact region, where with no fracture denoted as $\Omega_R:=\Omega \backslash \Omega_F\subset\Omega \in \mathbb{R}^d$ such that $\Omega_R\cup\Omega_F=\Omega$ and $\Omega_R\cap\Omega_F=\varnothing$. It has to be noted, that $\Omega_F$, i.e. the  domain in which the smeared crack phase-field is approximated, and its boundary $\partial \Omega_F$ strongly depend on the choice of the phase-field regularization parameter,
i.e. $\ep>0$. 

We deal with a multi-field problem depending on the displacement field $\bm u:\Omega\rightarrow\mathbb{R}^d$ and the crack phase-field $\varphi:\Omega\rightarrow[0,1]$ for a given pressure $p:\Omega\rightarrow\mathbb{R}$ and a given temperature $\Theta:\Omega\rightarrow\mathbb{R}$. The limiting values of $\varphi$, namely, $\varphi=1$ and $\varphi=0$ represent the undamaged and fully broken material phases, respectively. Finally, we, denote the $(\cdot, \cdot)$ the standard $L^2$ scalar product in  the domain $\Omega$ and $\langle \cdot, \cdot \rangle$ the $L^2$ scalar product on a boundary part.

\subsection{The energy functional in the thermo-poroelastic medium}
We start with the following energy functional posed in the thermo-poroelastic medium:
\begin{align}\label{lab_1}
     \mathcal{E} ( {\bm u}) &= \frac{1}{2}(\bm \sigma, \bm \varepsilon(
    {\bm u}))_{\Omega} 
- \langle\bm{\tau}, {\bm u}\rangle_{\partial_N\Omega},
\end{align}
where $\sigma$ is the stress tensor specified below, $\bm \varepsilon({\bm
  u})$ is the linearized strain tensor and $\bm{\tau}$ traction forces on 
the boundaries $\partial_N\Omega := \Gamma_N \cup \partial\mathcal{C}$. For the specific definition of the thermo-poroelastic stress, we 
apply the constitutive expression derived in \cite{Coussy2004}:
\[
\bm \sigma = \bm \sigma_0 + \bm \sigma_\varepsilon - \alpha_B (p-p_0)\bm{I} - 3\alpha_{\Theta} K_d(\Theta-\Theta_0)\bm{I}.
\]
Here, $\bm \sigma_0$ is the initial stress, $\alpha_B$ is Biot's constant, $p$
is the 
pressure, $p_0$ is the initial pressure, $\bm I$ is the second order identity tensor,
$3\alpha_\Theta$ is the 
volumetric skeleton thermal dilation coefficient (or thermal expansion
coefficient). Then, $K_d:=\frac{2}{d}\mu + \lambda$ ($d \in \{ 2,3\}$ still refers to the dimension of problem), 
$\Theta$ is the temperature, and $\Theta_0$ is the initial temperature, and 
\[
\bm \sigma_{\bm{\varepsilon}} = 2\mu {\bm{\varepsilon}} + \lambda tr({\bm{\varepsilon}})\bm I, \quad {\bm \varepsilon}(\bm u)=\frac{1}{2}(\nabla \bm u + \nabla \bm u^T),
\]
where $\mu$ and $\lambda$ are positive material parameters.

Using these definitions and setting $\bm\sigma_0 = 0$ (for the convenience of the
reader), we can write the energy functional in the thermo-poroelastic domain:
\begin{align}\label{lab_2}
     \mathcal{E} ( {\bm u}) &= \frac{1}{2}(\bm \sigma, {\bm{\varepsilon}}(
    {\bm u}))_{\Omega} 
- \langle\bm{\tau}, {\bm u}\rangle_{\partial_N\Omega}  \nonumber\\
&= 
\frac{1}{2}(\bm \sigma_{\bm \varepsilon}, {\bm{\varepsilon}}(
    {\bm u}))_{\Omega} 
- \langle{\bm \tau}, {\bm u}\rangle_{\partial_N\Omega} 
- (\alpha_B (p-p_0), \nabla\cdot \bm u)
- (3\alpha_\Theta K_d(\Theta-\Theta_0), \nabla\cdot \bm u).
\end{align}
This functional represents the energy in the thermo-poroelastic domain only.
For the total energy, we need to account for the crack energy as well.
In the presence of a fracture $\mathcal{C}$, we need to add 
the fracture energy which is expressed through \cite{FraMar98}:
\[
G_c \mathcal{H}^{d-1}  (\mathcal{C}),
\]
where $G_c$ is the critical elastic energy restitution rate
(related to the stress intensity factor through Irvins formula)
and $\mathcal{H}^{d-1}$ is the $d-1$-dimensional 
Hausdorff measure. 

Then, we arrive at the total energy functional:
\begin{align}\label{lab_3}
     \mathcal{E} ( {\bm u},\mathcal{C}) &= 
\frac{1}{2}(\bm \sigma_{\bm \varepsilon}, {\bm{\varepsilon}}({\bm u}))_{\Omega} 
- \langle{\bm \tau}, {\bm u}\rangle_{\partial_N\Omega}
 - (\alpha_B (p-p_0), \nabla\cdot \bm u)
- (3\alpha_\Theta K_d(\Theta-\Theta_0), \nabla\cdot \bm u) \nonumber \\
&\quad + G_c \mathcal{H}^{d-1}  (\mathcal{C}).
\end{align}
This functional is composed by contributions of two disjunct domains, namely 
$\Omega$ and $\mathcal{C}$. For the numerical treatment we regularize 
\eqref{lab_3} following \cite{BourFraMar00}. Specifically, 
the crack energy is approximated through a sequence of elliptic 
problems, so-called Ambrosio-Tortorelli functionals \cite{AmTo90,AmTo92}.
Therein, $\mathcal{H}^{d-1}$ is regularized by introducing 
an additional smoothed indicator variable nowadays called phase-field:
$\varphi:\Omega\to [0,1]$.
Finally, we account for the crack irreversibility constraint that the crack can only 
grow:
 \begin{equation}\label{evol1}
    \partial_t \varphi \leq 0.
 \end{equation}
For stating the variational formulations, we 
now introduce three sets:
\[
V := H^1_0 (\Omega)^d, \quad W:= H^1(\Omega), \quad W_{in} := \{ \varphi \in H^1(\Omega) | \; 0 \leq
\varphi\leq \varphi^{old} \}.
\]
As typical in problems with inequality constraints (see e.g.,
\cite{KiOd88,KiStam00}), $W_{in}$ is a nonempty, closed, convex,
subset of the linear function space $W$. Due to the inequality constraint \eqref{evol1}, $W_{in}$ is
not anymore a linear space.

Moreover, we recall that we do not deal with a classical time-dependent
problem, but a quasi-static load-incremental evolution, for which 
the loading interval is discretized using the discrete time (loading) points
\[
0 = t_0 < t_1 < \ldots < t_n < \ldots < t_N = T.
\]
With these preparations, the energy functional \eqref{lab_3} can be written as:

\newpage
\begin{form}[Regularized energy functional]
\label{form_1}
Let $p_0,p,\Theta,\Theta_0$ to be given with the initial conditions $\bm u_0=\bm u(x,0)$ and  $\varphi_0=\varphi(x,0)$. For the loading increments 
$n=1,2,\ldots, N$, find ${\bm u}:={\bm u}^n\in V$ and $\varphi:=\varphi^n\in W_{in}$ such that:
\begin{align}\label{lab_4}
    \mathcal{E}_{\ep} ( {\bm u}, \varphi) &= 
\underbrace{\frac{1}{2}(((1-\kappa)\varphi^2 + \kappa)\bm \sigma_{\varepsilon}, {\bm{\varepsilon}}(
    {\bm u}))_{\Omega}}_{\text{mechanical term}} 
- \underbrace{\langle\bm{\tau}, {\bm u}\rangle_{\Gamma_n \cup \partial\mathcal{C}}}_{\text{external term}} \nonumber \\ 
&\quad - \underbrace{(\alpha_B (p-p_0)\varphi^2, \nabla\cdot \bm u)}_{\text{pressure term}}
- \underbrace{(3\alpha_\Theta K_d(\Theta-\Theta_0)\varphi^2, \nabla\cdot \bm u)}_{\text{thermal term}} \nonumber \\
&\quad 
+ \underbrace{G_c \bigg( \frac{1}{2\ep} \|1-\varphi\|^2 + \frac{\ep}{2} \| \nabla \varphi
\|^2  \bigg)}_{\text{fracture term}}.
\end{align}
\end{form}

\begin{Remark}
Approximation of \eqref{lab_4} without pressures $p,p_0$ and temperatures
$\Theta, \Theta_0$ has been used in
many studies for solid mechanics.
\end{Remark}
\begin{Remark}
In \eqref{lab_4}, $\kappa := \kappa(h)$ 
is a positive regularization parameter for the elastic energy and 
$\ep :=\ep(h)$ is a regularization parameter for the phase-field 
variable denoting the width of the transition zone in which 
$\varphi$ changes from $0$ to $1$; see Fig. \ref{Figure1}.
\end{Remark}

The stationary points of the energy functional \eqref{lab_4} are characterized
by the first-order necessary conditions, namely the so-called Euler-Lagrange 
equations, which are obtained by (formal) differentiation
in the variables $\bm u$ and $\varphi$:

\begin{form}[Euler-Lagrange equations]
\label{form_EL_No_1}
Let $p_0,p,\Theta,\Theta_0$ be given with the initial conditions $\bm u_0=\bm u(x,0)$ and  $\varphi_0=\varphi(x,0)$. For the loading increments 
$n=1,2,\ldots, N$, find ${\bm u}:={\bm u}^n\in V$ and $\varphi:=\varphi^n\in W_{in}$ such that:
\begin{align*}
&(((1-\kappa)\varphi^2 + \kappa) \bm \sigma_{\bm \varepsilon}, {\bm{\varepsilon}}({\bm w})) 
- \langle {\bm \tau}, {\bm w}\rangle\\
&\quad - (\alpha_B (p-p_0)\varphi^2, \nabla \cdot {\bm w})
- (3\alpha_\Theta K_d(\Theta-\Theta_0)\varphi^2, \nabla\cdot {\bm w})
= 0 \quad\forall {\bm w}\in V,
\end{align*}
and
\begin{align*}
&(1-\kappa)(\varphi \bm \sigma_{\bm \varepsilon}(\bm u) : {\bm{\varepsilon}}(\bm u), \psi-\varphi)\\
&\quad - 2(\alpha_B -1) ((p-p_0) \nabla\cdot \bm u\, \varphi , \psi-\varphi)
- 2(3\alpha_\Theta K_d)((\Theta-\Theta_0)\nabla\cdot \bm u\, \varphi , \psi-\varphi)\\
&\quad +G_c (\frac{1}{\ep}(\varphi-1, \psi-\varphi) + \ep (\nabla\varphi,\nabla(\psi-\varphi))) \geq 0 \quad\forall \psi\in W \cap L^{\infty}.
\end{align*}
\end{form}

\subsection{Incorporating fracture pressure and temperature}
The energy functional \eqref{lab_4} and the corresponding
Euler-Lagrange equations in Formulation \ref{form_EL_No_1} 
are still incomplete because only mechanisms in the 
thermo-poroelastic medium are taken into account so far. In the following, we
consider pressure and temperature variations in the fracture 
as well, acting as additional normal stresses 
on the interface between the thermo-poroelastic medium and the fracture. 
Hence, we obtain a modification of the normal stress term 
$\langle {\bm \tau}, {\bm u}\rangle_{\partial_N \Omega}$.

Similar to \cite{MiWheWi18} in which an interface law 
was derived to account for pressurized fractures, we now
derive an expression including 
temperature effects. For the convenience of understanding, we also recapitulate 
the steps of the pressure and develop the modified terms for the temperature
simultaneously. Assuming that the fracture is a zone 
of high permeability and its width much smaller than its length, the 
leading order of the pressure stress in $\mathcal{C}$ is 
\begin{align*}
-(p-p_0){\bm I}.
\end{align*}
Similarly, we assume (according to \cite{TrSeNg13}) that the leading order 
of the temperature stress in $\mathcal{C}$ is given by 
\[
- C_\Theta (\Theta - \Theta_0){\bm I}, \quad\text{with } \; \; C_{\Theta} := C_{\Theta}(x,t).
\]
The constant $C_\Theta$ is obtained by following \cite{TrSeNg13} to get
a relation 
for the temperature near the fracture. We briefly sketch the basic idea.
Working with the heat conduction equation, while 
neglecting convection-dominated terms, we obtain:
\begin{equation}
\label{lab_5}
\partial_t \Theta - \nabla\cdot (\kappa_{\Theta} \nabla\Theta) = 0,
\end{equation}
with the thermal diffusivity
$\kappa_{\Theta} = \frac{K_r}{\rho_r C_r}$, where $\rho_r$ is the rock 
density, $C_r$ is the specific heat, $K_r$ thermal conductivity.
With the help of the heat conduction equation, 
a parameter measuring the 
gradient of the temperature at the 
fracture interface is derived.
This law provides an heuristic argument for 
the temperature evolution without
solving the full temperature equation in the whole domain.
Specifically, 
\begin{equation}
\label{C_T_eq}
C_\Theta = A_\Theta \left( \frac{\lambda_\Theta}{2\lambda_\Theta +1} \right),
\quad A_\Theta = \frac{E_Y\beta}{1-\nu_s},
\end{equation}
where $\beta$ is a linear thermal expansion coefficient, $E_Y$ is a Young's modulus,
$\nu_s$ is a Poisson's ratio. The parameter $\lambda_\Theta$ is derived 
by Hagoort \cite{Ha81} based on the heat conduction 
equation. Then,
\begin{equation}
\label{L_T_eq}
\lambda_\Theta = \sinh^{-1} (\frac{\gamma_T}{0.5l_0} \sqrt{\pi\kappa_{\Theta} t}),
\end{equation}
where $1\leq \gamma_T \leq \frac{4}{\pi}$.


\subsection{Interface laws for pressure and temperature}
We pursue in the following by deriving an interface relationship between
the thermo-poroelastic medium $\Omega$ and the fracture $\mathcal{C}$.
As usually required for coupled problems with interfaces, 
we need a kinematic and a dynamic coupling condition.
The kinematic conditions are to enforce 
\[
p_f = p_R \quad\text{and}\quad \Theta_f = \Theta_R \quad\text{on } \mathcal{C},
\]
where $p_f$ is the fracture pressure, $p_R$ the thermo-poroelastic pressure, 
$\Theta_f$ the temperature in the fracture and $\Theta_R$ the thermo-poroelastic temperature.
Since we enforce continuity of these variables on $\mathcal{C}$ we do not distinguish 
anymore in the rest of the paper and simply use $p$ and $\Theta$.

The dynamic coupling condition represents the continuity 
of normal stresses on the crack boundary $\partial\mathcal{C}$ and can be written as
\begin{equation}
\begin{aligned}
\label{eq_inter}
\bm \sigma {\bm n} &= [\underbrace{\bm \sigma_0 + \bm \sigma_{\bm \varepsilon} 
- \alpha_B (p-p_0) {\bm I} - 3\alpha_\Theta
K_d(\Theta-\Theta_0){\bm I}}_{\text{Thermo-poroelastic stress}}]{\bm n}\\
&= -\underbrace{[(p - p_0)- C_\Theta(\Theta - \Theta_0)]}_{\text{Fracture stress}}{\bm n}.
\end{aligned}
\end{equation}
Since the fracture boundary is smeared and not explicitly known, it remains 
to discuss how the interface condition is imposed. The kinematic 
conditions are Dirichlet-like conditions and build into the function spaces 
as usually done. The dynamic conditions are re-written as domain integrals by 
using Gauss' divergence theorem. We recall that the entire boundary 
is composed as
\[
\langle \bm\sigma\bm n, {\bm u}\rangle_{\partial_N \Omega}
 = \langle {\bm \tau}, {\bm u}\rangle_{\Gamma_N} + \langle \bm\sigma\bm n, {\bm u}\rangle_{\partial\mathcal{C}}.
\] 
We now manipulate the last term with the help of Eq. \ref{eq_inter} and the divergence theorem, reads
\begin{equation}
\begin{aligned}
\langle \bm \sigma {\bm n}, \bm u \rangle_{\partial\mathcal{C}} &= 
- \langle (p-p_0){\bm n}, \bm u\rangle_{\partial\mathcal{C}} 
- \langle C_\Theta(\Theta-\Theta_0)\bm{n}, \bm u\rangle_{\partial\mathcal{C}} \\
&= -(\nabla\cdot (p-p_0)\bm u) + \langle (p-p_0){\bm n}, \bm u\rangle_{\Gamma_N}
- (C_\Theta\nabla\cdot (\Theta-\Theta_0)\bm u) \\
& \quad + \langle C_\Theta(\Theta-\Theta_0){\bm n}, \bm u\rangle_{\Gamma_N}\\
&= - ((p-p_0),\nabla\cdot \bm u) - (\nabla (p-p_0),\bm u) + \langle (p-p_0)\bm{n},
\bm u\rangle_{\Gamma_N}\\
&\quad - (C_\Theta(\Theta-\Theta_0),\nabla\cdot \bm u) - (C_\Theta\nabla (\Theta-\Theta_0),\bm u) + \langle C_\Theta(\Theta-\Theta_0)\bm{n},
\bm u\rangle_{\Gamma_N}.
\end{aligned}
\end{equation}
These terms can be combined with the other terms containing pressures and temperatures in the energy
functional \eqref{lab_4}:
\[
(\alpha_B (p-p_0), \nabla\cdot \bm u) \quad\text{and}\quad 
(3\alpha_\Theta K_d(\Theta-\Theta_0), \nabla\cdot \bm u).
\]
In total we obtain (noticing that all signs change because we deal with 
$- \langle \bm\sigma\bm n, {\bm u}\rangle_{\partial_N \Omega}$ in \eqref{lab_4}):
\begin{equation}
\begin{aligned}
\label{eq_inter1}
&\quad - \langle \bm\sigma\bm n, {\bm u}\rangle_{\partial_N \Omega} - (\alpha_B (p-p_0),
\nabla\cdot \bm u) - (3\alpha_\Theta K_d(\Theta-\Theta_0), \nabla\cdot \bm u)\\
&= 
- \langle {\bm \tau},\bm u\rangle_{\Gamma_N} - \langle \bm \sigma n,\bm u\rangle_{\mathcal{C}} - (\alpha_B (p-p_0),
\nabla\cdot \bm u) - (3\alpha_\Theta K_d(\Theta-\Theta_0), \nabla\cdot \bm u)\\
&= - \langle {\bm \tau},\bm u\rangle_{\Gamma_N} - (\alpha_B (p-p_0), \nabla\cdot \bm u) +
((p-p_0),\nabla\cdot \bm u) + (\nabla (p-p_0),\bm u)\\
& \quad + \langle (p-p_0){\bm n},\bm u\rangle_{\Gamma_N} - (3\alpha_\Theta K_d(\Theta-\Theta_0), \nabla\cdot \bm u) + (C_\Theta(\Theta-\Theta_0),\nabla\cdot \bm u)\\ 
& \quad + (C_\Theta\nabla (\Theta-\Theta_0),\bm u) + \langle C_\Theta(\Theta-\Theta_0)\bm{n},\bm u\rangle_{\Gamma_N}.
\end{aligned}
\end{equation}
To complete, we introduce a phase-field representation of Eq. \ref{eq_inter1} and weight 
the domain terms with $\varphi^2$ (see Section 2 in \cite{MiWheWi18} for
pressurized fractures only)
in \eqref{lab_4}, which then yields:
\begin{form}[Energy functional including fracture pressure and temperature]
	\label{form_2}
Let $p_0,p,\Theta,\Theta_0$ be given with an initial condition $\bm u_0=\bm u(x,0)$ and  $\varphi_0=\varphi(x,0)$. For the loading increments 
$n=1,2,\ldots, N$, find $\bm u:=\bm u^n\in V$ and $\varphi:=\varphi^n\in W_{in}$ such that:
	\begin{align*}
	\mathcal{E}_{\ep} ( {\bm u},\varphi) &= 
    \underbrace{\frac{1}{2}(((1-\kappa)\varphi^2 + \kappa)\bm \sigma_{\varepsilon}(\bm u), {\bm \varepsilon}(
	{\bm u}))_{\Omega}}_{\text{Mechanical term}}
    - \underbrace{\bm{\langle\tilde\tau}, {\bm u}\rangle_{\Gamma_N}}_{\text{External load}} \\
	&\quad - \underbrace{((\alpha_B -1)(p-p_0)\varphi^2, \nabla\cdot {\bm u}) + (\nabla
	(p-p_0)\varphi^2,{\bm u})}_{\text{Pressure term}}\\
	&\quad - \underbrace{((3\alpha_\Theta K_d + C_\Theta)(\Theta - \Theta_0)\varphi^2, \nabla\cdot {\bm u})
	+ (C_\Theta \nabla (\Theta-\Theta_0)\varphi^2, {\bm u})}_{\text{Thermal term}}\\
	&\quad + \underbrace{G_c \bigg( \frac{1}{2\ep} \|1-\varphi\|^2 + \frac{\ep}{2} \| \nabla \varphi
	\|^2  \bigg)}_{\text{Fracture term}},
	\end{align*}
	where $ \tilde{\bm \tau} = {\bm \tau} - (p-p_0){\bm n} + C_\Theta(\Theta-\Theta_0){\bm n}$. 
\end{form}

\begin{Remark}
	We underline again that $p$ and $\Theta$ are fixed and given in this paper
	and therefore not solution variables. The extension in which
	$p$ is also an unknown was first proposed in \cite{MiWheWi14}.
        A similar development for non-isothermal phase-field fractures in
        thermo-poroelasticity is ongoing work.
\end{Remark}
\begin{Remark}
	We notice that the temperature opens the crack if $\Theta < \Theta_0$ (see also Remark \ref{Rem1}), i.e., the current (or injected)
	temperature is lower than the initial temperature.
\end{Remark}
\begin{Remark}
When normal stresses (traction forces) are prescribed on (parts of) the
boundary, the term $ \tilde{\bm \tau} = {\bm \tau} - (p-p_0){\bm n} +
C_\Theta(\Theta-\Theta_0){\bm n}$ must be carefully considered; see
\cite{SogoLeeWheeler_2018}
or \cite{MiWheWi18}[Section 5.2].
\end{Remark}

\newpage
\section{Thermodynamic arguments and mechanical analysis}
\label{sec_thermo}
In this section, we further augment our proposed model with thermodynamic arguments and subsequently analyze our model 
from a mechanical point of view.

\subsection{Extension to a decoupled strain-energy function into volumetric and isochoric response} \label{sec_decoupled_energy}
Since the fracturing material behaves quite differently in bulk and shear
parts of the domain, we employ a consistent split for the strain energy
density function, i.e. 
\[
\Psi(\bm\varepsilon(\bm u)):=\frac{\lambda}{2}\big(\bm\varepsilon(\bm u):\bm
I\big)^2+\mu {\bm \varepsilon}(\bm u)^2:{\bm I}.
\]
Hence, instead of dealing directly with $\bm\varepsilon(\bm u)$, we perform
additive decomposition of strain tensor into \textit{volume-changing}
(volumetric part) and \textit{volume-preserving} (deviatoric part), i.e.
\[
\bm\varepsilon(\bm u)=\bm\varepsilon^{vol}(\bm u)+\bm\varepsilon^{dev}(\bm
u).
\]
Here, the volumetric strain is denoted as $\bm\varepsilon^{vol}(\bm
u):=\frac{1}{3}(\bm \varepsilon(\bm u):\bm I)\bm I$ and the deviatoric strain is
denoted as $\bm\varepsilon^{dev}(\bm u):=\mathbb{P}:\bm \varepsilon$. The
fourth-order projection tensor
$\mathbb{P}:=\mathbb{I}-\frac{1}{3}\bm I\otimes\bm I$ 
is introduced to map the full strain tensor to its deviatoric counterpart.
Therein,
$\mathbb{I}_{i,j,k,l}:=\frac{1}{2}\big(\delta_{i,k}\delta_{j,l}+\delta_{i,l}\delta_{j,k}\big)$
is the fourth-order symmetric identity tensor. Furthermore, 
$\mathbb{P}$ possesses the major symmetries, i.e. $\mathbb{P}_{i,j,k,l}=\mathbb{P}_{k,l,i,j}$, and $\mathbb{P}^n=\mathbb{P}$ for any given integer $n$. So, a decoupled representation of the strain-energy function into a so-called volumetric and deviatoric contribution are given as follows,

\begin{equation}
{\Psi}(\bm\varepsilon(\bm u))={\Psi}(\bm\varepsilon^{vol}(\bm u))+{\Psi}(\bm\varepsilon^{dev}(\bm u)).
\label{Q1}
\end{equation}
Therein, the volumetric contribution of the strain energy density function reads,

\begin{equation}
{\Psi}(\bm\varepsilon^{vol}(\bm u)):=\frac{\lambda}{2}\big(\bm\varepsilon^{vol}(\bm u):\bm I\big)^2+\mu {\bm \varepsilon}^{vol}(\bm u)^2:{\bm I}=\frac{K_d}{2}(\bm\varepsilon^{vol}(\bm u):\bm I)^2,
\label{Q2}
\end{equation}
where $K_d:=\frac{2}{d}\mu + \lambda$ is the bulk modulus and $d\in \{2,3\}$. 
The deviatoric contribution of the strain energy density function is

\begin{equation}
{\Psi}(\bm\varepsilon^{dev}(\bm u)):=\frac{\lambda}{2}\big(\bm\varepsilon^{dev}(\bm u):\bm I\big)^2+\mu {\bm \varepsilon}^{dev}(\bm u)^2:{\bm I}=\mu {\bm \varepsilon}^{dev}(\bm u)^2:{\bm I}.
\label{Q3}
\end{equation}
To show that equality holds in Eq.\ref{Q1}, the identities $\bm\varepsilon^{dev}(\bm u):\bm I=0$ and $\bm\varepsilon^{vol}(\bm u):\bm\varepsilon^{dev}(\bm u)=0$ are used. Physically, it is trivial to assume that the degradation induced by the phase field acts only on the tensile and shear counterpart of the elastic strain density function. Hence, it is assumed there is no degradation in compression, which also prevents interpenetration of the crack lips during crack closure \cite{AmGeraLoren15}. It turns out that the modified strain energy density function for the fracturing material becomes,

\begin{equation}
\Psi\big(\bm\varepsilon(\bm u)\big):=g(\varphi){\Psi^+}\big(\bm\varepsilon(\bm u)\big)+{\Psi^-}\big(\bm\varepsilon(\bm u)\big),
\label{Q4}
\end{equation}
such that a monotonically decreasing quadrature degradation function, i.e. 
\begin{equation}
g(\varphi):=(1-\kappa)\varphi^2 + \kappa,
\label{N4}
\end{equation}
describes the degradation of the solid with the evolving crack phase-field
parameter $\varphi$, includes $\kappa$ that is a small residual stiffness that
is introduced to prevent numerical problems. 
Additionally, the positive part of the strain energy density function that is the tensile and deviatoric part of full strain energy density function reads

\begin{equation}
{\Psi^{+}}(\bm\varepsilon(\bm u))=H{^+}(\nabla. \bm u){\Psi}(\bm\varepsilon^{vol}(\bm u))+{\Psi}(\bm\varepsilon^{dev}(\bm u)).
\label{Q5}
\end{equation}
Therein, $H{^+}(\nabla. \bm u)$ is a \textit{positive Heaviside function} such that if $\nabla. \bm u$ is positive return one and otherwise, give zero value. It is noted due to identity $\nabla. \bm u=\text{tr}(\bm \varepsilon)$, the positive Heaviside function indicates the points in a domain where they are in tensile part, i.e. $H{^+}(\nabla. \bm u)=1$ and compression part, i.e. $H{^+}(\nabla. \bm u)=0$. Negative strain energy density function that is the compression part of the full strain energy density function, is

\begin{equation}
{\Psi^{-}}(\bm\varepsilon(\bm u))=\big(1-H{^+}(\nabla. \bm u)\big){\Psi}(\bm\varepsilon^{vol}(\bm u)).
\label{Q6}
\end{equation}
The constitutive equation for the modified strain energy density function, whereas the stress response constitutes an additive split of $\bm \sigma_{\bm \varepsilon}:=\frac{\partial \Psi(\bm \varepsilon)}{\partial \bm \varepsilon}$ to purely tensile contribution, i.e. ${\bm \sigma}^{+}_{\bm \varepsilon}({\bm{\varepsilon}})$ and a purely compression contribution, i.e. ${\bm \sigma}^{-}_{\bm \varepsilon}({\bm{\varepsilon}})$, reads

\begin{equation}
{\bm \sigma}_{\bm \varepsilon}:=\frac{\partial \Psi({\bm \varepsilon}}{\partial {\bm \varepsilon}}=g(\varphi)\frac{\partial \Psi^+({\bm \varepsilon})}{\partial \bm \varepsilon}+\frac{\partial \Psi^-({\bm \varepsilon})}{\partial \bm \varepsilon}=g(\varphi){\bm \sigma^{+}_{\bm \varepsilon}}({\bm{\varepsilon}})+{\bm \sigma^{-}_{\bm \varepsilon}}({\bm{\varepsilon}}).
\label{Q7}
\end{equation}
Therein,

\begin{align*}
&{\bm \sigma}^{+}_{\bm \varepsilon}({\bm{\varepsilon}})=K_nH{^+}(\nabla. \bm u)(\bm \varepsilon:\bm I)\bm I+2\mu{\bm \varepsilon}^{dev},\quad \text{and} \quad {\bm \sigma^{-}_{\bm \varepsilon}}({\bm{\varepsilon}})=K_n\big(1-H{^+}(\nabla. \bm u)\big)(\bm \varepsilon:\bm I)\bm I.
\label{Q8}
\end{align*}
The decoupled representation of the fourth-order elasticity tensor to relate the work into conjugate pairs of stress and strain tensor by means of additive decomposition of the stress tensor, reads

\begin{equation}
\mathbb{C}:=\frac{\partial {\bm \sigma}( \bm{\varepsilon})}{\partial {\bm \varepsilon}}=g(\varphi)\frac{\partial {\bm \sigma}^{+}( \bm{\varepsilon})}{\partial {\bm \varepsilon}}+\frac{\partial {\bm \sigma}^{-}( \bm{\varepsilon})}{\partial {\bm \varepsilon}}=:g(\varphi)\mathbb{C}^++\mathbb{C}^-,
\label{Q9}
\end{equation}
with

\begin{align*}
&\mathbb{C}^+=K_nH{^+}(\nabla. \bm u)\bm I\otimes\bm I+2\mu\mathbb{P},\quad \text{and} \quad \mathbb{C}^-=K_n\big(1-H{^+}(\nabla. \bm u)\big)\bm I\otimes\bm I,
\label{Q10}
\end{align*}
where the identity $\frac{\partial \bm \varepsilon^{dev}(\bm u)}{\partial \bm \varepsilon(\bm u)}=\mathbb{P}$ is used.
\begin{form}[Final energy functional with strain-energy density splitting]
	\label{form_3}
	Let $p_0,p,\Theta_0,\Theta$ be given with the initial conditions $\bm u_0=\bm u(x,0)$ and  $\varphi_0=\varphi(x,0)$. For the loading increments 
$n=1,2,\ldots, N$, find ${\bm u}:={\bm u}^n\in V$ and $\varphi:=\varphi^n\in W_{in}$ such that:
	\begin{align*}
	\mathcal{E}_{\ep} ( {\bm u},\varphi) &= 
	\underbrace{\frac{1}{2}\Big(g(\varphi_+)\bm {\sigma^+}_{\bm \varepsilon}(\bm u), {\bm \varepsilon}(
		{\bm u})\Big)_{\Omega}+(\bm {\sigma^-}_{\bm \varepsilon}(\bm u), {\bm \varepsilon}({\bm u}))_{\Omega}}_{\text{Mechanical term}}
	- \underbrace{\bm{\langle\tilde\tau}, {\bm u}\rangle_{\Gamma_N}}_{\text{External load}} \\
	&\quad - \underbrace{((\alpha_B -1)(p-p_0)\varphi_+^2, \nabla\cdot {\bm u}) + (\nabla
		(p-p_0)\varphi_+^2,{\bm u})}_{\text{Pressure term}}\\
	&\quad - \underbrace{((3\alpha_\Theta K_d + C_\Theta)(\Theta - \Theta_0)\varphi_+^2, \nabla\cdot {\bm u})
		+ (C_\Theta \nabla (\Theta-\Theta_0)\varphi_+^2, {\bm u})}_{\text{Thermal term}}\\
	&\quad + \underbrace{G_c \bigg( \frac{1}{2\ep} \|1-\varphi\|^2 + \frac{\ep}{2} \| \nabla \varphi
		\|^2  \bigg)}_{\text{Fracture term}}.
	\end{align*}
\end{form}

\begin{Remark}
	\label{Rem_phi_pos}
    In the case of elastic cracks, it can be shown that the phase field unknown satisfies $0 \leq \varphi \leq 1$. In order
    to establish this property for the spatially continuous incremental problem, we need to modify energy functional for the
    negative values of $\varphi$. Hence, similar to \cite{MiWheWi18}[Section 2],
   we have used $\varphi_+$ rather than $\varphi$ in terms where negative $\varphi$ could lead to incorrect
    physics in the bulk energy, traction, pressure and thermal forces.    
\end{Remark}

\begin{Remark}
	\label{alternariveE_AT1_AT2}
We briefly mention 	
	alternative descriptions for the phase-field approximation and the
        degradation function.
First, using a more general term regarding fracture term is denoted as
$\frac{G_c}{4c_w}( \frac{1}{2\ep} W(\varphi) + \frac{\ep}{2} \| \nabla \varphi
\|^2 )$ for the constant $c_w$ and $W(\varphi)$, which refers to the local part
of the dissipated fracture energy functional. This is the so-called
$\operatorname{\texttt{AT-1}}$ fracture model if $c_w=\frac{2}{3}$ and
$W(\varphi)=\|1-\varphi\|$, see \cite{Pham_AT1}, and it is
$\operatorname{\texttt{AT-2}}$ for $c_w=\frac{1}{2}$ and
$W(\varphi)=\|1-\varphi\|^2$; see \cite{Bourdin_AT2}. 
Second, the quadratic polynomial $g(\varphi)$ can be written in the form of the cubic polynomial as 	$g(\varphi):=3(1-\kappa)\varphi^2-2\varphi^3 + \kappa$ or quartic polynomial as $g(\varphi):=4(1-\kappa)\varphi^3-3\varphi^4 + \kappa$, see \cite{kuhn_g_c}.
\end{Remark}


For $\mathcal{E}_{\ep} ( {\bm u},\varphi)$ given in Formulation \ref{form_3}, we derive the characterizing Euler-Lagrange equations by
differentiating the energy functional with respect to ${\bm u}$ and $\varphi$:
\begin{form}[Final Euler-Lagrange equations]
	\label{form_4}
	Let $p_0,p,\Theta_0,\Theta$ be given with the initial conditions $\bm u_0=\bm u(x,0)$ and  $\varphi_0=\varphi(x,0)$. 
	For $n=1,2,3,\ldots, N$, find $\bm u:= \bm u^n\in V$:
	\begin{align*}
	A_1 ( {{\bm u}})({\bm w}) &= (g(\varphi_+){\bm \sigma^+_{\bm
            \varepsilon}}(\bm u), {\bm \varepsilon}(\bm w)) + ({\bm \sigma^-_{\bm \varepsilon}(\bm u)}, {\bm \varepsilon}(\bm w)) 
	- \langle \tilde {\bm \tau}, {\bm w}\rangle_{\Gamma_N}\\
	&\quad - (\alpha_B - 1)((p-p_0)\varphi_+^2, \nabla\cdot {\bm w}) 
	+ (\nabla (p-p_0)\varphi_+^2, {\bm w})\\
	&\quad - (3\alpha_\Theta K_d + C_\Theta)((\Theta-\Theta_0)\varphi_+^2 , \nabla\cdot {\bm w})
	+ (C_\Theta \nabla (\Theta_F - \Theta_0)\varphi_+^2, {\bm w})
	= 0 \quad\forall {\bm w}\in V,
	\end{align*}
	and find $\varphi:=\varphi^n\in W_{in}$:
	\begin{align*}
	A_2 (\varphi)(\psi-\varphi) &=(1-\kappa)(\varphi_+ {\bm \sigma^+_{\bm
            \varepsilon}}(\bm u) : {\bm \varepsilon}(\bm u), \psi-\varphi)\\
	&\quad - 2(\alpha_B -1) (\varphi_+(p-p_0) \nabla\cdot {\bm u}\,  , \psi-\varphi)
	+ 2(\varphi_+ \nabla (p-p_0){\bm u}\, ,\psi-\varphi)\\
	&\quad - 2(3\alpha_\Theta K_d + C_\Theta)(\varphi_+(\Theta-\Theta_0)\nabla\cdot {\bm u}\, , \psi-\varphi)
	+ 2(C_\Theta \varphi_+ \nabla (\Theta-\Theta_0) {\bm u}\, , \psi-\varphi)\\
	&\quad + G_c \left( \frac{1}{\ep}(\varphi-1, \psi-\varphi)
	+ \ep (\nabla\varphi,\nabla(\psi-\varphi))\right) \geq 0 \quad\forall
        \psi\in W \cap L^{\infty}.
	\end{align*}

\end{form}


\subsection{The Euler-Lagrange equations in a strong form}
\label{sec_strong_form}

In order to complete our derivations, we derive the strong form of Formulation
\ref{form1} in this section. Using integration by parts, we obtain a quasi-stationary elliptic system for the displacements and the phase-field variable, where the latter one is subject to an inequality constraint in time and therefore needs to be complemented with a complementary condition:

\begin{form}[Strong form of the Euler-Lagrange equations]
	\label{form1}
Let a pressure $p: \Omega \rightarrow \mathbb{R}$ and temperature $\Theta: \Omega
\rightarrow \mathbb{R}$ and the initial conditions $\bm u_0=\bm u(x,0)$ and
$\varphi_0=\varphi(x,0)$ be given. 
For the loading steps $n=1,2,3,\ldots,
N$, we solve a displacement equation where we seek $\bm u:= \bm u^n: \Omega \rightarrow \mathbb{R}^{d}$
	\begin{equation}
	\begin{aligned}
	&-\nabla. \big( g(\varphi_+)\bm {\sigma^+}_{\bm \varepsilon}(\bm u)+\bm {\sigma^-}_{\bm \varepsilon}(\bm u)\big)\\
	&+(\alpha_B - 1)\nabla.(\varphi_+^2(p-p_0))+\varphi_+^2\nabla(p-p_0)\\
	&+(3\alpha_\Theta K_d + C_\Theta)\nabla.(\varphi_+^2(\Theta-\Theta_0))
+\varphi_+^2\nabla(\Theta-\Theta_0)= 0\quad in \; \Omega
        \end{aligned}
        \end{equation}
\begin{equation}
 {\bm u} =\bm 0 \quad on \; \partial\Omega.
	\end{equation}
	The phase-field system consists of three parts: the PDE, the
        inequality constraint, and a compatibility condition (in fracture
        mechanics called Rice condition \cite{ricecond}).
    Find $\varphi:=\varphi^n : \Omega \rightarrow [0,1]$
    \begin{equation}
\label{strong_2}
    \begin{aligned}
    &-\Big( ((1-\kappa)\varphi_+ \Psi^+(\bm\varepsilon(\bm u))
    -G_c\ep\Delta\varphi
    -\frac{G_c}{\ep}(1-\varphi)\\
    &- 2(\alpha_B -1) \varphi_+(p-p_0) \nabla\cdot {\bm u}+ 2\varphi_+ \nabla (p-p_0){\bm u}\\
    &- 2(3\alpha_\Theta K_d + C_\Theta)\varphi_+(\Theta-\Theta_0)\nabla\cdot {\bm u}
    + 2C_\Theta \varphi_+ \nabla (\Theta-\Theta_0) {\bm u}
    \Big)\leqslant 0 \quad in \; \Omega\\	
    \end{aligned}
    \end{equation}

\begin{equation}
\label{eq_ine_con}
  \partial_t\varphi\leqslant 0 \quad in \; \Omega,
\end{equation}

	\begin{equation}
\label{strong_4}
    \begin{aligned}
    &-\Big( ((1-\kappa)\varphi_+ \Psi^+(\bm\varepsilon(\bm u))
    -G_c\ep\Delta\varphi_+
    -\frac{G_c}{\ep}(1-\varphi_+)\\
    &- 2(\alpha_B -1) \varphi_+(p-p_0) \nabla\cdot {\bm u}+ 2\varphi_+ \nabla (p-p_0){\bm u}\\
    &- 2(3\alpha_\Theta K_d + C_\Theta)\varphi_+(\Theta-\Theta_0)\nabla\cdot {\bm u}
    + 2C_\Theta \varphi_+ \nabla (\Theta-\Theta_0) {\bm u}
    \Big) \; \partial_t\varphi= 0\quad in \; \Omega,\\
    \end{aligned}
    \end{equation}
    
    \begin{equation}
    {\partial_{\bm n} {\varphi}} = 0 \quad on \; \partial\Omega.
    \end{equation}
    
\end{form}


\subsection{Global balance principle of continuum thermo-mechanics}
In this section, thermodynamical consistency for the preservation of the principal of balance of energy is shown by considering a sequence of variational substitutions. As a point of departure, considering Formulation \ref{form_3}, we define the total free energy functional as

\begin{equation}
\begin{aligned}
\mathcal{E}(\bm u,\varphi_+) &:= \mathcal{E}_{bulk}(\bm u,\varphi_+)+\mathcal{E}_{crack}(\varphi)+\mathcal{E}_{ext}(\bm u)\\
&:=\int_\Omega W_{bulk}(\bm\varepsilon(\bm u),\varphi_+) \, \mathrm{d}{\textbf{x}}+\int_\Omega W_{crack}(\varphi) \, \mathrm{d}{\textbf{x}}           
-\int_{\Gamma_{N} } {\bm {\tilde\tau}} \cdot \bm u\,\mathrm{d}s,
\label{N1}
\end{aligned}
\end{equation}
with the bulk free energy functional

\begin{equation}
\begin{aligned}
W_{bulk}\bigg(\bm\varepsilon(\bm u),\varphi_+\bigg):=
& g(\varphi_+){\Psi^+}(\bm\varepsilon(\bm u))+{\Psi^-}(\bm\varepsilon(\bm u))\\
& - ((\alpha_B -1)(p-p_0)\varphi_+^2. \nabla\cdot \bm u) + (\nabla
(p-p_0)\varphi_+^2.\bm u)\\
& - ((3\alpha_\Theta K_d + C_\Theta)(\Theta - \Theta_0)\varphi_+^2. \nabla\cdot \bm u)
+ (C_\Theta \nabla (\Theta-\Theta_0)\varphi_+^2. \bm u),\\
\label{N2}
\end{aligned}
\end{equation}
and the crack free energy functional
\begin{equation}
W_{crack}(\varphi,\nabla \varphi):= G_c \gamma_{\ep}(\varphi,\nabla \varphi).
\label{N3}
\end{equation}
Here, we have introduced the second order crack surface density function per unit volume of the thermo-poroelastic media as,

\begin{equation}
\gamma_{\ep}(\varphi,\nabla \varphi):=  \frac{1}{2\ep} |1-\varphi|^2 + \frac{\ep}{2} | \nabla \varphi
|^2 .
\label{N6}
\end{equation}
We recall that weak forms derived for the displacement and phase-field are given in Formulation \ref{form_4}. To derive a global balance of energy, we choose as test functions $\bm w=:\dot {\bm u}$ and $\psi=:\dot {\varphi}$ and restate Formulation \ref{form_4} for the mechanical part as,
\begin{equation}
\begin{aligned}
A_1 ( {{\bm u}})(\bm {\dot u})=&(g(\varphi_+)\bm \sigma^+_{\bm \varepsilon}({\bm u}), \nabla(\bm{\dot u}))+(\bm \sigma^-_{\bm \varepsilon}({\bm u}), \nabla(\bm{\dot u})) 
- \langle \tilde{\bm \tau}, \bm{\dot u}\rangle_{\Gamma_N}\\
&\quad - (\alpha_B - 1)((p-p_0)\varphi_+^2, tr(\nabla\bm{\dot u}) 
+ (\nabla (p-p_0)\varphi_+^2,\bm{\dot u})\\
&\quad - (3\alpha_\Theta K_d + C_\Theta)((\Theta-\Theta_0)\varphi_+^2 , tr(\nabla\bm{\dot u}))
+ (C_\Theta \nabla (\Theta - \Theta_0)\varphi_+^2, \bm{\dot u})
= 0 \quad\forall \bm{\dot u}\in V.
\label{N7}
\end{aligned}
\end{equation}
We define the second order displacement gradient denoted as $\bm H:=\nabla{\bm u}$, then Eq. (\ref{N6}) reduces to 

\begin{equation}
\begin{aligned}
A_1 ( {{\bm u}})(\bm {\dot u})=&(g(\varphi_+)\bm \sigma^+_{\bm \varepsilon}(\bm u), \bm{\dot H}) +(\bm \sigma^-_{\bm \varepsilon}(\bm u), \bm{\dot H})
- \langle \tilde{\bm \tau}, \bm{\dot u}\rangle_{\Gamma_N}\\
&\quad - (\alpha_B - 1)((p-p_0)\varphi_+^2, tr(\bm{\dot H}) 
+ (\nabla (p-p_0)\varphi_+^2,\bm{\dot u})\\
&\quad - (3\alpha_\Theta K_d + C_\Theta)((\Theta-\Theta_0)\varphi_+^2 , tr(\bm{\dot H}))
+ (C_\Theta \nabla (\Theta - \Theta_0)\varphi_+^2, \bm{\dot u})
= 0 \quad\forall \bm{\dot u}\in V,
\label{N8}
\end{aligned}
\end{equation}
that is

\begin{equation}
\begin{aligned}
\int_\Omega \dot W_{bulk}\bigg(\bm u,\bm\varepsilon(\bm u),\varphi_+\bigg)\, \mathrm{d}{\textbf{x}}-\int_\Omega \frac{\partial W_{bulk}}{\partial \varphi} \dot \varphi \, \mathrm{d}{\textbf{x}}-P^{ext}=0
\label{N9},
\end{aligned}
\end{equation}

\begin{equation}
\begin{aligned}
\Leftrightarrow \qquad P^{int}-\int_\Omega \frac{\partial W_{bulk}}{\partial \varphi} \dot \varphi \, \mathrm{d}{\textbf{x}}-P^{ext}=0.
\label{N10}
\end{aligned}
\end{equation}
Herein, the rate of internal mechanical power which describes the response of a domain $\Omega$ done by the stress field is 

\begin{equation}
P^{int}=\int_\Omega \dot W_{bulk}\bigg(\bm u,\bm\varepsilon(\bm u),\varphi_+\bigg)=\int_\Omega \frac{\partial W_{bulk}}{\partial \bm { u}}. \bm {\dot u} \, \mathrm{d}{\textbf{x}}
+\int_\Omega \frac{\partial W_{bulk}}{\partial \bm { \varepsilon}} : \bm{\dot \varepsilon} \, \mathrm{d}{\textbf{x}}
+\int_\Omega \frac{\partial W_{bulk}}{\partial \varphi} \dot \varphi \, \mathrm{d}{\textbf{x}}.
\label{N11}
\end{equation}
Accordingly, through Formulation \ref{form_4} for the phase-field part, we derive

\begin{equation}
A_2 (\varphi)(\dot \varphi)=\int_\Omega \frac{\partial W_{bulk}}{\partial \varphi} \dot \varphi \, \mathrm{d}{\textbf{x}}+ \int_\Omega G_c \left( \frac{1}{\ep}(\varphi-1). \dot \varphi+ \ep (\nabla\varphi)\nabla\dot \varphi\right) = 0 \quad \quad \forall \dot \varphi>0\in W.
\label{N12}
\end{equation}
It is important to note that the inequality $A_2(\cdot)(\cdot)$ in Formulation
\ref{form_4} 
becomes an equality in Eq. \ref{N12} because we consider the situation in which the 
inequality constraint  \eqref{eq_ine_con} is strictly fulfilled; namely
$\dot \varphi < 0$. In this case equality must hold in \eqref{strong_2} (and
thus in $A_2(\cdot)(\cdot)$) since 
otherwise the compatibility condition \eqref{strong_4} is not fulfilled. 
We can restate Eq. \ref{N12} as follows,

\begin{equation}
\dot {\mathcal{E}}_{crack}(\varphi)= \int_\Omega G_c \delta_{\varphi} \gamma_{\ep}(\varphi,\nabla \varphi) \dot \varphi \mathrm{d}{\textbf{x}} = -\int_\Omega \frac{\partial W_{bulk}}{\partial \varphi} \dot \varphi \, \mathrm{d}{\textbf{x}},
\label{N13}
\end{equation}
where the left-hand side is considered to be a functional of the rate of the crack phase-field and that is the global crack dissipation functional, i.e. $\dot {\mathcal{E}}_{crack}(\varphi)$, hence

\begin{equation}
\dot {\mathcal{E}}_{crack}(\varphi) \geqslant 0.
\label{N14}
\end{equation}
Through the second term in Eq. \ref{N13} and considering Eq. \ref{N14}, the local form of crack dissipation reads,

\begin{equation}
\delta_{\varphi} \gamma_{\ep}(\varphi,\nabla \varphi) \leqslant 0 \qquad \text{and} \qquad \dot \varphi \leqslant 0,
\label{N141}
\end{equation}
and the third term of Eq. \ref{N13} leads to the additional local condition, i.e. 

\begin{equation}
\beta_{\varphi}:=\frac{\partial W_{bulk}(\bm u, \varphi_+)}{\partial \varphi} \geqslant 0.
\label{N19}
\end{equation}
Herein, $\beta_{\varphi}$ introduced as a crack deriving force is conjugate to the phase-field variable, that is

\begin{equation}
\begin{aligned}
\beta_{\varphi}&=2(1-\kappa)\varphi_+ \Psi(\bm\varepsilon(\bm u))\\
& - 2((\alpha_B -1)(p-p_0)\varphi_+, \nabla\cdot \bm u) + 2(\nabla
(p-p_0)\varphi_+,\bm u)\\
& - 2((3\alpha_\Theta K_d + C_\Theta)(\Theta - \Theta_0)\varphi_+, \nabla\cdot \bm u)
+ 2(C_\Theta \nabla (\Theta-\Theta_0)\varphi_+, \bm u).\\
\label{N20}
\end{aligned}
\end{equation}
In our formulation, in comparison to Miehe et al. \cite{MieWelHof10b} we have $\beta_{\varphi}=-f$ (due to the definition of the crack phase-field like the damage variable in the Miehe et al. \cite{MieWelHof10b}). Note that $\beta_{\varphi}$ becomes zero in the fracture surface, i.e. $\mathcal{C}$, as $\varphi$ becomes zero. 

\begin{Remark}
	\label{Rem1}
     {We now look into more detail in the positivity condition for the
       Eq. \ref{N20} and discuss term by term.
It is noted due to the condition obtained from Eq. \ref{N14} which is the
positivity of the global crack dissipation functional, and by considering
third term in Eq. \ref{N13} leads to the positivity of the crack deriving
force, i.e. $\beta_{\varphi} \geqslant 0$  and $\dot \varphi \leqslant 0$,
i.e. Eq. \ref{N19}. Hence, for the evolving fracture, by means of
Eq. \ref{N20}, the first term including $\Psi(\bm\varepsilon(\bm u))$ is
positive due to the normalization condition and for the second and fourth
terms we are  demanding for pressure $p \geqslant p_0$ and thermal $\Theta
\leqslant \Theta_0$ both conditions hold (i.e. $\beta_{\varphi} \geqslant 0$)
to satisfy positiveness of the global crack dissipation functional (it is
assumed the gradient term for both pressure and thermal part are neglected and
$\alpha_B\ll 1$). It has to be noted, due to decoupled strain energy density
function we have used (see Section \ref{sec_decoupled_energy}),  results in
$\nabla. \bm u >  0$ in $\Omega_F$ (i.e., the transition zone). This can be
observed later in Fig. \ref{Figure13} and \ref{Figure15}. Additionally,
$\dot \varphi \leqslant 0$ is imposed to the
total energy functional through the primal-dual active set strategy (see
Section \ref{sec_primal_dual_as}). 
By doing so, the first
condition in Eq. \ref{N141} is automatically satisfied (due to satisfaction of
the third term in Eq. \ref{N13}. This results in the positivity
of the second term in Eq. \ref{N13})} if and only if both $p \geqslant p_0$ and $\Theta \leqslant \Theta_0$ hold.
\end{Remark}
To derive the Karush-Kuhn-Tucker conditions, 
we consider Eq. \ref{N13}, and the functional derivative with respect to the $\varphi$ defined as:

\begin{equation}
\int_\Omega \delta_{\varphi} \gamma_{\ep}(\varphi,\nabla \varphi) \mathrm{d}{\textbf{x}}:=\int_\Omega \bigg( \frac{\partial \gamma}{\partial \varphi} +\frac{\partial \gamma}{\partial \nabla \varphi}\bigg) \mathrm{d}{\textbf{x}} =
\int_\Omega \bigg( \frac{\partial \gamma}{\partial \varphi}-\nabla.[\frac{\partial \gamma}{\partial \nabla \varphi}]\bigg) \mathrm{d}{\textbf{x}}, 
\label{N15}
\end{equation}
that is

\begin{equation}
\int_\Omega \delta_{\varphi} \gamma_{\ep}(\varphi,\nabla \varphi) \mathrm{d}{\textbf{x}}=\int_\Omega \frac{1}{\ep}[(\varphi-1)-\ep^2 \Delta \varphi] \mathrm{d}{\textbf{x}}.
\label{N16}
\end{equation}
Due to $\frac{\partial \gamma_{\ep}(\varphi,\nabla \varphi)}{\partial \nabla
  \varphi}.\bm n=\ep \nabla \varphi.\bm n=0$, that is Euler-Lagrange crack
phase-field equation on the boundary of the given domain. The right hand side
of Eq. \ref{N13}, i.e. $-\int_\Omega \beta_{\varphi} \dot \varphi \,
\mathrm{d}{\textbf{x}}$, can also be obtained by means of Eq. $\ref{N10}$
and hence Eq. \ref{N13} is restated as,

\begin{equation}
\Pi(\dot {\bm u}, \dot {\varphi})= \dot {\mathcal{E}}_{crack}(\varphi)+P^{int}(\dot {\bm u}, \dot {\varphi})-P^{ext}(\dot {\bm u}) = 0,
\label{N17}
\end{equation}
that is a global form of the balance of energy for the coupled two-field problem describing the evolution of internal energy and crack dissipation energy, i.e. $\dot {\mathcal{E}}_{crack}(\varphi)+P^{int}$, in a system due to the external loads, i.e. $P^{ext}$. Additionally, by means of the second and the third terms in Eq. \ref{N13}, it turns out that  

\begin{equation}
\int_\Omega (\beta_{\varphi}+G_c \delta_{\varphi} \gamma_{\ep}(\varphi,\nabla \varphi)) \dot \varphi \mathrm{d}{\textbf{x}} = 0,
\label{N18}
\end{equation}
which is the balance law for the evolution of the crack phase-field which
ensuring the principal of maximum dissipation during the crack phase-field
evolution (see e.g. \cite{MieWelHof10b}) and so-called as a compatibility
condition; we refer the reader also to Section \ref{sec_strong_form}.
Observe that one may satisfy this global irreversibility constraint of crack evolution, i.e. Eq. \ref{N18} by ensuring locally a positive variational derivative of the crack surface function and a positive evolution of the crack phase field, i.e.

\begin{equation}
\beta_{\varphi}+G_c \delta_{\varphi} \gamma_{\ep}(\varphi,\nabla \varphi) \leqslant 0, \quad \text{and} \quad \dot \varphi \leqslant 0.
\label{N21}
\end{equation}

The former condition is ensured in the subsequent treatment by a constitutive assumption that relates the functional derivative to a positive energetic driving force. The latter constraint is a natural assumption that relates the fracture phase field for the non-reversible evolution of crack phase-field.

\begin{Remark}
	\label{compatibility_const}
It is noted within loading state, i.e. $\dot \varphi < 0$, due to the compatibility condition, i.e. Eq. \ref{N18}, along with Eq. \ref{N21}, one may observe  $\beta_{\varphi}+G_c \delta_{\varphi} \gamma_{\ep}(\varphi,\nabla \varphi) = 0$ and in the unloading state, i.e. $\dot \varphi = 0$, we have $\beta_{\varphi}+G_c \delta_{\varphi} \gamma_{\ep}(\varphi,\nabla \varphi) < 0$. Equation \ref{N21} along with Eq. \ref{N18} refer to the Karush-Kuhn-Tucker conditions for the phase-field fracturing problem. 
\end{Remark}

\begin{Remark}
	\label{Rem2}
	It turns out that, if the positivity of the crack deriving force is satisfied (see Remark \ref{Rem1}), i.e. $\beta_{\varphi} \geqslant 0$, and according to Eq. \ref{N21} we have $G_c \delta_{\varphi} \gamma_{\ep}(\varphi,\nabla \varphi) \leqslant  -\beta_{\varphi} $  and hence results to $\delta_{\varphi} \gamma_{\ep}(\varphi,\nabla \varphi) \leqslant 0$, that is the first inequality condition shown in Eq. \ref{N141}.
\end{Remark}

\section{Numerical solution and the final discrete model}
\label{sec_numerics}
In this section, we briefly describe spatial discretization first. 
The solution algorithm is then based on a quasi-monolithic approach for which 
a Newton solver is employed as described in \cite{HeWheWi15}. 
The crack irreversibility condition in Eq. \ref{evol1}
is treated with a primal-dual active set method. Both techniques can 
be gathered into one single combined Newton
solver \cite{HeWheWi15,LeeWheWi16,HeiWi18}.

\subsection{Spatial discretization}
The computational domain is subdivided into 
quadrilateral or hexahedral element domains. Both subproblems are discretized with a Galerkin finite element method 
using $H^1$-conforming bilinear (2D) or trilinear (3D) elements, i.e., the
ansatz and test space uses $Q_1^c$-finite elements, e.g., for details, we refer readers to the \cite{Cia87}.
Consequently, the discrete spaces have the property 
$V_h \subset V$ and $W_h \subset W$.

\subsection{Problem statement of the compact minimization problem}
\label{sec_setup}
The starting point for the discrete solution is Formulation \ref{form_3},
written 
in compact form as:
\begin{eqnarray*}
&&\min_{\bm u\in V, \varphi\in W_{in}} \mathcal{E}_\varepsilon (\bm u,\varphi), 
\end{eqnarray*}
For the following, we set $U = (\bm u, \varphi) \in V \times W$. Discretizing
\[
 \partial_t \varphi \approx \frac{\varphi^{n+1} - \varphi^n}{\delta t},
\]
with the time step size $\delta t:=t^{n+1}-t^n$,
the incremental problem can be rewritten as
\begin{equation}
\label{cons_min_Energy}
\begin{aligned}
&\min \mathcal{E}_{\ep} ({\bm U}) \\
&\textup{subject to }  {\bm U} \leq \bm {U}^{old} \textup{ on } \Phi,
\end{aligned}
\end{equation}
where $\Phi = 0 \times W$, so that the constraint acts on the phase-field
variable only, and $\bm {U}^{old}$ is the solution from the last time-step (or the
initial condition).
The minimization of \eqref{cons_min_Energy} is numerically challenging, due to the following reasons:

\begin{enumerate}
	\item the energy functional $\mathcal{E}_{\ep} ( {\bm u},\varphi)$ 
may admit several local minimizers. Thus finding the global minimum is in
general non-feasible, e.g., \cite{BourFraMar08} for discussions 
on the pure elasticity case. The existence of a minimizer 
for pressurized fractures was established in \cite{MiWheWi18}. Since 
the temperature is a given quantity, as the pressure, the proof for existence
of the current $\mathcal{E}_{\ep} ( {\bm u},\varphi)$ goes along the same
lines as in \cite{BourFraMar08},
	\item the irreversibility of crack phase-field, i.e. $|\mathcal{C}_{t-1}|\leq |\mathcal{C}_{t}|$, is required to provide a thermodynamically consistent minimization problem by having a positive crack dissipation inequality and  enforcing on the temporal derivative of the phase-field function, see e.g. \cite{MieWelHof10},
	\item the minimization problem is characterized by localization of the
          crack phase-field in bands of width  of order $\ep$. From a
          practical and numerical analysis point of view, $\eps > h$ must
          hold (at least one element has to be existed to cover regularized
          phased-field). In more detail, $h=o(\eps)$.  The regularization parameter is
          typically a very small dimensionless value and for the accurate
          fracture response (i.e. converge toward the sharp crack profile)
          which should tend to $0$ in the limit $h \rightarrow 0$ 
          to resolve the bands, see e.g. \cite{NoiiGL},
	\item the linear system of equations arises from Hessian matrix of the $\mathcal{E}_{\ep} ( {\bm u},\varphi)$ are typically badly conditioned due to the presence of crack phase-field localizations band where the elastic stiffness varies rapidly from the intact value to zero, see e.g. \cite{FarrellSVPFF}.
	
\end{enumerate}


\subsection{A combined Newton method: treating crack irreversibility and
  solving the nonlinear problem}
\label{sec_primal_dual_as}
In following, we will address how to resolve the issues mentioned above. To
this end, we first describe Newton's method for solving the unconstrained minimization problem $\min \mathcal{E}_\ep(\bm U)$ in Eq. \ref{cons_min_Energy} for the total energy functional $ \mathcal{E}_\ep(\bm U)$ given in formulation \ref{form_3}.
We construct a sequence ${\bm U}^0, {\bm U}^1, \dots, {\bm U}^N$ with
\[
 {\bm U}^{k+1} = {\bm U}^k + \delta {\bm U}^{k},
\]
where the update $\delta {\bm U}^k$ is computed as the solution of the linear system
(details in Section \ref{sub_sec_monolithic_form}):
\begin{align}
  \label{eq:outer-newton}
  \nabla^2 \mathcal{E}_\ep({\bm U}^k) \, \delta {\bm U}^{k} = -\nabla \mathcal{E}_\ep ({\bm U}^k).
\end{align}
If we assume the constraints on the phase-field on Eq. \ref{cons_min_Energy} hold for the initial 
guess ${\bm U}^0$ (we will start with the solution from the last time step, which satisfies the constraint), the condition
\begin{align}
 \label{eq:outer-constraint}
 \delta {\bm U}^k \leq \textbf{0} \textup{ on } \Phi, 
\end{align}
implies that ${\bm U}^{k+1} = {\bm U}^k + \delta {\bm U}^{k} \leq {\bm U}^k \leq \dots 
\leq {\bm U}^0 \leq {\bm U}^{old} \textup{ on } \Phi$.
In a variational formulation, the previous Newton method reads:
\begin{equation}
\begin{aligned}
  \label{linear_system_in_Newton}
   \nabla^2 \mathcal{E}_\ep( {{\bm U}_h^k})(\delta {\bm U}_h^k, {\bm \Psi})  = - \nabla \mathcal{E}_\ep( {{\bm U}_h^k})( {\bm \Psi}), \quad\forall{\bm \Psi}\in V_h\times W_h, \quad \text{with } \quad \delta {\bm U}^k \leq 0 \textup{ on } \Phi,
   \end{aligned}
\end{equation}
where ${\bm \Psi}:=[\bm w, \psi]$ is denoted as the total test function. 
Crack irreversibility is taken care of by removing the corresponding 
rows and columns in which the constraint is is active. This yields then 
a reduced system.
In our implementation, we combine two Newton methods (active set and the nonlinear iteration for the PDE solution; see Section \ref{sub_sec_monolithic_form}) into a single update loop with variable $\delta {\bm U}^k$. This Newton loop contains a back-tracking line search to improve the convergence radius.
This yields Algorithm \ref{final_algo}.

\newpage
\begin{table}[!ht]\small
	\caption{\em Combined Newton loop at time step $t^n$.}
	\label{final_algo}
	{
		\begin{tabular}{l}
			\hline \\
			{\textbf{Input}:} loading data $(\bm u, \varphi, p_{t},p_0,\Theta,\Theta_0,t)$ on $\Gamma_{D},\Gamma_{N}\subset\partial\Omega$; \\ 
			\hspace{1.25cm}solution ${\bm U}_{n-1}:=(\bm u_{n-1},\varphi_{n-1})$ from time step $n-1$. \\ [0.1cm]
			\quad\quad Combined Newton / primal-dual active set iteration $k\geq1$:\\
			\quad\quad\quad \textbullet\; Assemble residual ${\bm R}({\bm U}_h^k)$,\\
			\quad\quad\quad \textbullet\; Compute active set $\mathcal{A}_k = \{i \mid  ({\bm B}^{-1})_{ii} ({\bm R}_k)_i + c (\delta {\bm U}_h^k)_i > 0 \}$,\\
			\quad\quad\quad \textbullet\; Assemble matrix ${\bm G} = \nabla^2 E_\ep({\bm U}_h^k)$ and right-hand side ${\bm F} = -\nabla E_\ep({\bm U}_h^k)$, \\
			\quad\quad\quad \textbullet\; Eliminate rows and columns in $\mathcal{A}_k$ from $\bm G$ and $\bm F$ to obtain $\widetilde{\bm G}$
			and $\widetilde{\bm F}$, \\
			\quad\quad\quad \textbullet\; Solve linear system $\widetilde{\bm G} \delta \bm U_k = \widetilde{\bm F}$, i.e, find $\delta \bm U_h^k \in V_h \times W_h$ with, \\
			$\quad \quad \quad \quad \quad \quad \; \: \nabla^2 \mathcal{E}_\ep( {{\bm U}_h^k})(\delta {\bm U}_h^k, {\bm \Psi})  = - \nabla \mathcal{E}_\ep( {{\bm U}_h^k})( {\bm \Psi}), \quad\forall{\bm \Psi}\in V_h\times W_h,$
			\\ \quad \quad \quad \; \; where $\nabla^2 \mathcal{E}_\ep$ and $\nabla \mathcal{E}_\ep$ are defined in Section \ref{sub_sec_monolithic_form} and \ref{sec_block_structure}.\\
			\quad\quad\quad \textbullet\; Find a step size $0 < \omega \leq 1$ using back-tracking line search algorithm to get\\
			\quad \quad \quad \quad \quad \quad \quad \; ${\bm U}_h^{k+1} = {\bm U}_h^k + \omega \delta {\bm U}_h^k$
			with $ \widetilde{\bm R}({\bm U}_h^{k+1})<\widetilde{\bm R}({\bm U}_h^k)$. \\
			\quad\quad\quad \textbullet\; if fulfilled, such that\\
			\quad \quad \quad \quad \quad \quad \quad \; $ \mathcal{A}_{k+1} = \mathcal{A}_k \; \text{and} \; \widetilde{\bm R}({\bm U}_h^k) <\operatorname{\texttt{TOL}_\texttt{N-R}}\widetilde{\bm R}({\bm U}_h^0$)\\
			\quad \quad \quad \; \; set $(\bm u^k,\varphi^k)=:(\bm u_t,\varphi_t)$ and stop; \\
			\quad\quad\quad {\color{white}\textbullet}\; else $k+1\rightarrow k$. \\ [0.1cm]
			{\textbf{ Output:}} solution $(\bm u_n,\varphi_n)=:{\bm U}_{n}$. \\ \\
			\hline
		\end{tabular}
	}
\end{table}


\begin{Remark}[Stopping criteria]
	\label{Stopping_criteria}
We remark that the above algorithm has two stopping criteria that must be 
achieved simultaneously:
\begin{equation}
\mathcal{A}_{k+1} = \mathcal{A}_k \qquad \text{and} \qquad
\widetilde{\bm R}({\bm U}_h^k) < \operatorname{\texttt{TOL}_\texttt{N-R}} \widetilde{\bm R}({\bm U}_h^0).
\end{equation}
\end{Remark}

\begin{Remark}
 It is important to distinguish between the full residual ${\bm R}({\bm U}^k_h)$ and $\widetilde{\bm R}({\bm U}_h^k)$.
 The latter is the residual on the inactive set, which can be computed by eliminating the active set constraints from the former.
\end{Remark}


\subsection{On the Jacobian and the residual inside Newton's method}
\label{sub_sec_monolithic_form}

For  solving  the  PDE  problem at each Newton step,
we  focus  on  a  monolithic  scheme  in  which  all  equations  are  solved simultaneously resulting in one semi-linear form. However, it is well known that the energy functional (3) is not convex simultaneously in both solution variables $\bm u$ and $\varphi$; but separately in each variable while keeping the other fixed. Consequently, solving the Euler-Lagrange equations in a straightforward way is not possible and influences
the robustness and efficiency of the solution scheme because of an (possibly)
indefinite Hessian matrix $G$ \cite{GeLo16,Wi17_SISC,Wi17_CMAME}.

\begin{figure}
	\centering
	{\includegraphics[clip,trim=4cm 6cm 4cm 1cm, width=15cm] {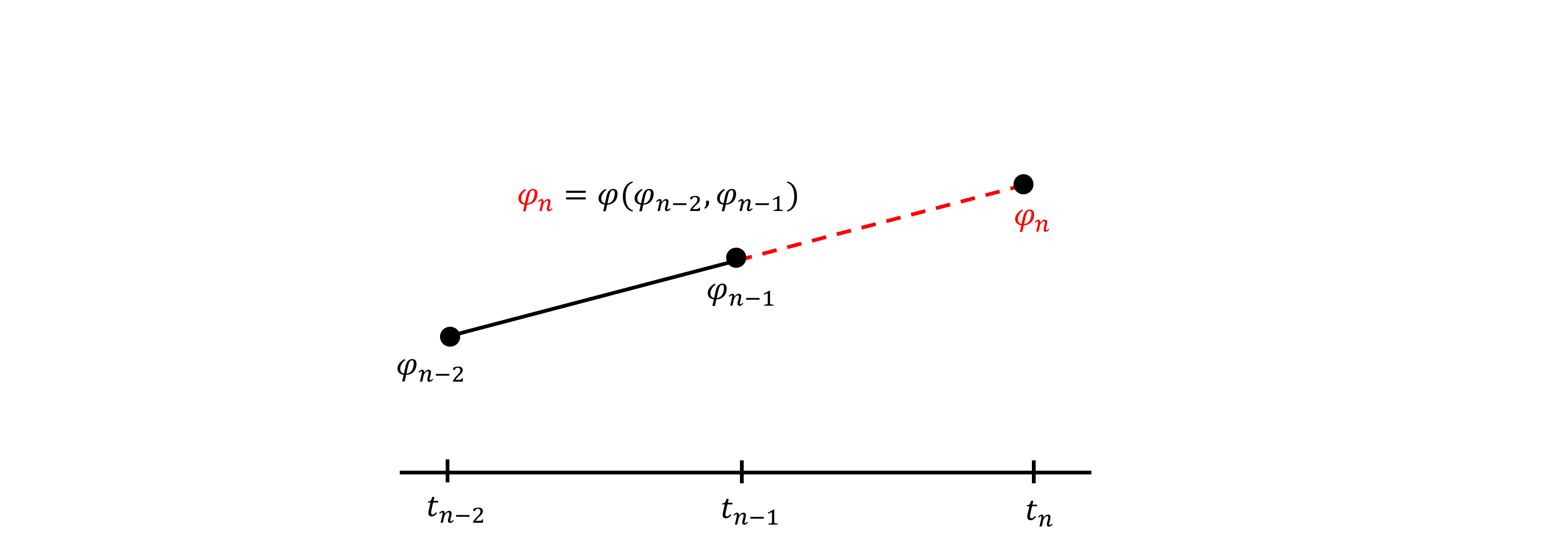}}   
	\caption{Linear extrapolation in time for the $\varphi$ based on two previous solutions, i.e. $\varphi^{n-2},\varphi^{n-1}$ .}
	\label{Figure2}
\end{figure}

The critical terms arise in the Hessian matrix $G$, are the cross terms,
includes for mechanical term i.e. 
$\big(((1-\kappa){ \varphi_+^2} + \kappa)\bm \sigma^+_{\bm \varepsilon}(\bm u), {\bm \varepsilon}(\bm w)\big)$,
for pressure term i.e. $(\alpha_B - 1)((p-p_0){\varphi_+^2}, \nabla\cdot \bm w)$,
$(\nabla (p-p_0) \varphi_+^2, \bm w)$ and for thermal term i.e. $(3\alpha_\Theta K + C_\Theta)((\Theta-\Theta_0){ \varphi_+^2}, \nabla\cdot \bm w)$ and
$(C_\Theta \nabla (\Theta - \Theta_0){ \varphi}^2, \bm w)$. To this end, according to \cite{HeWheWi15}, for having a convex energy functional, $A_1$ is linearized  in the direction of the  $\bm u$ and $\varphi$ by linear extrapolation and time-lagging of the phase-field, i.e. $\varphi \approx \tilde \varphi:= \varphi(\varphi^{n-2},\varphi^{n-1})$ in order to obtain a convex energy functional, that is
\begin{equation}
\tilde \varphi:= \varphi(\varphi^{n-2},\varphi^{n-1})=\varphi^{n-2}\frac{t^n-t^{n-1}}{t^{n-2}-t^{n-1}}+\varphi^{n-1}\frac{t^n-t^{n-2}}{t^{n-1}-t^{n-2}}.
\end{equation}
Here, $\varphi^{n-2},\varphi^{n-1}$ denote the solutions to previous time
steps, denoted as $t^{n-2}$ 
and $t^{n-1}$, see Fig.~\ref{Figure2}.

Hence, linearization of  $A_1(\bm u)(\bm w)$ (see Formulation \ref{form_4}) in the direction of $\varphi$, i.e. $\delta \varphi$, within tangent stiffness matrix (the Hessian matrix) is neglected (due to the given $\tilde \varphi$).

In the following, we state monolithic formulations for the
displacement-phase-field system. The presentation is similar 
to \cite{HeWheWi15}. Specifically, the phase-field variable 
is time-lagged in the first term of the displacement equation.
For fully monolithic solution algorithms and their performance we refer 
the reader to \cite{Wi17_SISC,Wi17_CMAME,GeLo16}.
In the following, we deal with equalities since we removed the 
inequality constraint via the primal-dual active set strategy.
The residual reads:
\begin{equation}
\label{monolithic_semi_linear}
\begin{aligned}
A({\bm U})( {\bm \Psi}) := { \nabla \mathcal{E}_\ep({\bm U})}( {\bm \Psi})
&=(g(\tilde \varphi_+)\bm \sigma^+_{\bm \varepsilon}(\bm u), {\bm \varepsilon}(\bm w))+(\bm \sigma^-_{\bm \varepsilon}(\bm u), {\bm \varepsilon}(\bm w)) 
- \bm{\langle \tilde\tau}, \bm w\rangle_{\Gamma_N}\\
\quad &- (\alpha_B - 1)((p-p_0){\tilde \varphi_+}^2, \nabla\cdot \bm w) 
+ (\nabla (p-p_0){\tilde \varphi_+}^2, \bm w)\\
\quad &- (3\alpha_\Theta K_d + C_\Theta)((\Theta-\Theta_0){\tilde \varphi_+}^2 , \nabla\cdot \bm w)
+ (C_\Theta \nabla (\Theta - \Theta_0){\tilde \varphi_+}^2, \bm w)\\
\quad&+(1-\kappa)(\varphi_+ \bm \sigma^+_{\bm \varepsilon}(\bm u) : {\bm \varepsilon}(\bm u), \psi)\\
\quad&- 2(\alpha_B -1) ((p-p_0) \nabla\cdot \bm u\, \varphi_+ , \psi)
+ 2(\nabla (p-p_0)\bm u\, \varphi_+ ,\psi)\\
\quad&- 2(3\alpha_\Theta K_d + C_\Theta)((\Theta-\Theta_0)\nabla\cdot \bm u\, \varphi_+, \psi)
+ 2(C_\Theta \nabla (\Theta-\Theta_0) \bm u\, \varphi_+ , \psi)\\
\quad&+ G_c \left( \frac{1}{\ep}(\varphi-1, \psi)
+ \ep (\nabla\varphi,\nabla\psi)\right).
\end{aligned}
\end{equation}
The corresponding Jacobian is built by computing 
the directional derivative $A'( {\bm U})(\delta {\bm U}, {\Psi})$. Then:
\begin{equation}
\label{monolithic_semi_linear_prime}
\begin{aligned}
A'(\bm U) (\delta {\bm U}, {\Psi}) :&= \nabla^2 \mathcal{E}_\ep(\bm U) (\delta {\bm U}, {\Psi}) \\ 
&=\Bigl( 
g(\tilde \varphi_+)\;\bm \sigma^+_{\bm \varepsilon}({ \delta \bm u}), {\bm \varepsilon}( \bm w )\Bigr)+(\sigma^-_{\bm \varepsilon}({ \delta \bm u}), {\bm \varepsilon}( \bm w ))\\
&+ (1-\kappa) \big(\delta\varphi_+ \bm \sigma^+(\bm u):{\bm \varepsilon}(\bm u) + {\varphi_+} \;
\bm \sigma^+_{\bm \varepsilon}({ \delta \bm u}):{\bm \varepsilon}( {\bm u})+{\varphi_+} \;
\bm \sigma^+_{\bm \varepsilon}({ \bm u}):{\bm \varepsilon}( \bm \delta u),\psi\big) \\
&-  2(\alpha_B -1) (p-p_0) (\delta\varphi_+\mbox{div } \bm u
+ {\varphi_+}\; \mbox{div }  { \delta \bm u},\psi) \\
&+  2\nabla (p-p_0) (\delta\varphi \bm u
+ {\varphi_+}\; { \delta \bm u},\psi) \\
&-  2(3\alpha_\Theta K_d + C_\Theta)(\Theta-\Theta_0) (\delta\varphi\mbox{div } \bm u
+ {\varphi_+}\; \mbox{div }  { \delta \bm u},\psi) \\
&+  2C_\Theta\nabla (\Theta-\Theta_0) (\delta\varphi \bm u
+ {\varphi_+}\; { \delta \bm u},\psi) \\
&+  G_c  \Bigl( \frac{1}{\ep} (\delta\varphi,\psi) + \ep (\nabla
\delta\varphi, \nabla \psi)   \Bigr).
\end{aligned}
\end{equation}

\newpage
\subsection{On the linear equation system at each Newton step}
\label{sec_block_structure}
\subsubsection{Spatial discretization and block structure}

In this section, we consider the structure and solution of the 
linear discrete system \eqref{linear_system_in_Newton} arising in each Newton step. 
For spatial discretization, we use the previously introduced 
spaces $V_h\times W_h$ with vector-valued basis
\[
\{ \psi_i \, | i = 1,\ldots,N_s\},
\]
where the basis functions are primitive (they are only non-zero in one component), so we can
separate them into displacement and phase-field basis functions and sort them accordingly:
\begin{eqnarray*}
 \psi_i &= 
 \begin{pmatrix}  \chi^{\bm u}_i \\ 0 \end{pmatrix}
 , \text{ for } i=1,\dots,N_u, \\
 \psi_{(N_{\bm u}+i)} &=
 \begin{pmatrix} 0 \\ \chi^\varphi_i \end{pmatrix},
 \text{ for } i=1,\dots,N_\varphi,
\end{eqnarray*}
where $N_{\bm u} + N_\varphi = N_s$. 
This is now used to 
transform \eqref{linear_system_in_Newton} into a system of the
form 
\begin{equation}
\label{lin_sys_Ax_b}
{\bm M}{\bm x} = {\bm F},
\end{equation}
where ${\bm M}$ is a block matrix (the Jacobian) and 
$F$ the right-hand side consisting of the residuals.
The block structure is
\[
 {\bm M} = \begin{pmatrix}
      {\bm M}^{{\bm u} {\bm u}} & {\bm M}^{{\bm u}\varphi} \\
      {\bm M}^{\varphi {\bm u}} & {\bm M}^{\varphi\varphi} \\
     \end{pmatrix},
     \qquad
 {\bm F} = \begin{pmatrix}
      F^{\bm u} \\ F^\varphi
     \end{pmatrix},
\]
with entries coming from \eqref{monolithic_semi_linear_prime}:
\begin{equation}
\begin{aligned}
 {\bm M}^{{\bm u}{\bm u}}_{i,j}
&=
\Bigl(\big( (1-\kappa) {\tilde\varphi_+}^2  +\kappa \big)\;{\bm \sigma}^+_{\bm \varepsilon}({ \chi}_j^{\bm u}), {\bm \varepsilon}( { \chi}_i^{\bm u} )\Bigr)+({\bm \sigma}^-_{\bm \varepsilon}({\chi}_j^{\bm u}), {\bm \varepsilon}( { \chi}_i^{\bm u} )),\\
{\bm M}^{\varphi {\bm u}}_{i,j}
&= 
2(1-\kappa) ({\varphi_+} \; {\bm \sigma}^+_{\bm \varepsilon}({   \chi}_j^{\bm u}):{\bm \varepsilon}( {\bm u}),\chi_i^{\varphi})\\
&-  2(\alpha_E -1) (p-p_0) ({\varphi_+}\; \nabla. ( { \chi}_j^{\bm u}),\chi_i^{\varphi})\\
&+  2 \nabla(p-p0) {\varphi_+}\; ( { \chi}_j^{\bm u},\chi_i^{\varphi})\\
&-  2(3\alpha_\Theta K_d + C_\Theta)((\Theta-\Theta_0) ({\varphi_+}\; \nabla. ( { \chi}_j^{\bm u}),\chi_i^{\varphi})  \\
&+  2C_\Theta \nabla (\Theta-\Theta_0) {\varphi_+}\; ( { \chi}_j^{\bm u},\chi_i^{\varphi}),\\
{\bm M}^{{\bm u}\varphi}_{i,j}
&=0, \\
{\bm M}^{\varphi\varphi}_{i,j}
&= (1-\kappa) ( {\bm \sigma}^+_{\bm \varepsilon}(\bm u):{\bm \varepsilon}(\bm u) \chi_j^{\varphi},\chi_i^{\varphi})\\
&\quad -  2(\alpha_B -1) (p-p_0) (\nabla.( \bm u) \chi_j^{\varphi},\chi_i^{\varphi}) \\
&\quad  -  2(3\alpha_\Theta K_d + C_\Theta)((\Theta-\Theta_0) (\nabla.( \bm u) \chi_j^{\varphi},\chi_i^{\varphi}) \\
&\quad +  2{\bm u} \nabla (p-p_0) (\chi_j^{\varphi},\chi_i^{\varphi}) \\
&\quad +  2{\bm u} C_\Theta \nabla (\Theta-\Theta_0) (\chi_j^{\varphi},\chi_i^{\varphi}) \\
&\quad+  G_c  \Bigl( \frac{1}{\ep} (\chi_j^{\varphi},\chi_i^{\varphi}) + \ep (\nabla
\chi_j^{\varphi}, \nabla \chi_i^{\varphi})   \Bigr). 
\end{aligned}
\end{equation}

\begin{Remark}
	Since we replaced $\varphi^2$ by $\tilde\varphi^2$ in the displacement equation,
	the block $M^{{\bm u}\varphi}_{i,j}$ is zero
	and the Hessian matrix $\bm G$ has triangular structure. In the other case, all blocks 
	would be nonzero; see \cite{Wi17_SISC,Wi17_CMAME}.
\end{Remark}

The right-hand side consists of the corresponding 
residuals (see semi-linear form \eqref{monolithic_semi_linear}). 

In particular, we have
\begin{equation}
\begin{aligned}
{\bm F}^{\bm u}_i &= - \widetilde{A}({\bm U}_k)(\chi_i^{\bm u})\\ 
&= \Bigl(\big( (1-\kappa) {\tilde\varphi_{+k}}^2  +\kappa \big)\;{\bm \sigma}^+_{\bm \varepsilon}({\bm u}_k), {\bm \varepsilon}( {{\chi_i^{\bm u}}} )\Bigr)+({\bm \sigma}^-_{\bm \varepsilon}({\bm u}_k), {\bm \varepsilon}( {{\chi_i^{\bm u}}} ))  
- \langle \tilde {\bm \tau}, {{{\chi_i^{\bm u}}}}\rangle_{\Gamma_N}\\ 
&\;\;\;\;- (\alpha_B - 1) ({\tilde\varphi}_{+k}^{2} (p-p_0), {\nabla. }  {\chi_i^{\bm u}})
+(\nabla (p-p_0){\tilde \varphi}_{+k}^{2}, {\chi_i^{\bm u}})\\
&\;\;\;\;- (3\alpha_\Theta K_d + C_\Theta) ({\tilde\varphi}_{+k}^{2} (\Theta-\Theta_0), {\nabla. }  {\chi_i^{\bm u}})
+C_\Theta(\nabla (\Theta-\Theta_0){\tilde \varphi}_{+k}^{2}, {\chi_i^{\bm u}}),\\\\
{\bm F}^\varphi_i &= - \widetilde{A}({\bm U}_k)(\chi_i^{\varphi}) \\
&= (1-\kappa) ({\varphi}_{+k} \;{\bm \sigma}^+_{\bm \varepsilon}({\bm u}_k):\bm{\varepsilon}( {\bm u}_k),
\chi_i^\varphi) \\
&\;\;\;\; -  2(\alpha_B -1) ({\varphi_k}\;  (p-p_0)\; {\nabla. }  {\bm u}_k,\chi_i^\varphi) 
+ 2({\varphi_{+k}}\;  \nabla (p-p_0)\; {\bm u}_{+k},\chi_i^\varphi)\\
&\;\;\;\; -  2(3\alpha_\Theta K_d + C_\Theta) ({\varphi_{+k}}\;  (\Theta-\Theta_0)\; {\nabla. }  {\bm u}_k,\chi_i^\varphi) 
+ 2(C_\Theta{\varphi_{+k}}\;  \nabla (\Theta-\Theta_0)\; {\bm u}_k,\chi_i^\varphi)\\
&\;\;\;\; +  G_c  \Bigl( -\frac{1}{\ep} (1-\varphi_k,\chi_i^\varphi) + \ep (\nabla
\varphi_k, \nabla \chi_i^\varphi)   \Bigr).
\end{aligned}
\end{equation}

The criterion for convergence of the Newton method is based on the relative residual norm that is, $\texttt{Residual}: \| \bm F(\bm x_{k+1}) \| \leq \texttt{Tol}_\texttt{N-R} \| \bm F(\bm x_{k}) \|$ for the user-prescribed $\texttt{Tol}_\texttt{N-R}$.

In the matrix, the degrees of freedom that belong to Dirichlet conditions
(here only displacements since we assume Neumann conditions for the phase-field)
are strongly enforced by replacing the corresponding rows and columns as
usual in a finite element code. In a similar fashion, the rows 
and columns that belong to nodes of the active set are removed from the
matrix. Corresponding right-hand side values in the vector $F$ are set to
zero; see as well Table  \ref{final_algo}.

\subsubsection{Solution of the linear equation systems}
\label{sec_lin_sol}
The linear system arising at each Newton step is solved iteratively using a generalized minimal residual (GMRES) scheme with a block diagonal preconditioner, i.e. ${\bm P}^{-1}$, as follows

\[
\begin{pmatrix}
{\bm M}^{{\bm u}{\bm u}} & \bm 0 \; \; \\
{\bm M}^{\varphi {\bm u}} & {\bm M}^{\varphi\varphi} \\
\end{pmatrix}
\begin{pmatrix}
F^{\bm u} \\ F^\varphi
\end{pmatrix}=
\begin{pmatrix}
	\delta {\bm u} \\ \delta \varphi
\end{pmatrix} \quad \text{with} \quad {\bm P}^{-1}:=
\begin{pmatrix}
	\big({ \tilde{\bm M}}^{{\bm u}{\bm u}}\big)^{-1} & {\bm 0} \\
	{\bm 0} & \big(\tilde{\bm M}^{\varphi \varphi}\big)^{-1} \\
\end{pmatrix},
\]	
where we approximate the blocks using a single V-cycle of algebraic multigrid (by Trilinos ML, \cite{trilinos1}), 
and due to the in-dependency to the mesh element size, the whole solver scheme is nearly optimal as recently demonstrated for pressurized phase-field fractures in \cite{HeiWi18}. The main reason why this solver works nicely is that the block-diagonal terms 
are of elliptic type.

\section{Numerical studies}
\label{sec_tests}
The performance of the proposed model is investigated by means of a 
series of representative numerical examples in two- and three spatial dimensions.  
The proposed formulation is considered to be a canonically consistent and robust scheme for treating
pressurized fractures and non-isothermal setting in thermo-poroelastic media. 
The emphasis is on verifications of the programming code against 
analytical solutions obtained with Sneddon-Lowengrub's formula
\cite{SneddLow69} (extended to non-isothermal configurations 
in \cite{TrSeNg13}) and using realistic material properties.
Mesh refinement studies are provided to study the accuracy of our approach.
Investigations of the effects strain-energy splitting and solver outcomes are
studied in addition.
The implementation is based on deal.II \cite{dealII85,BangerthHartmannKanschat2007} and 
specifically on the programming code from \cite{HeiWi18}.

\subsection{Geometries and parameters}
\label{sec_geometry_param}

\subsubsection{Two-dimensional setup}
\begin{figure}
	\centering
	{\includegraphics[trim=0cm 3cm 0cm 3cm, width=15cm] {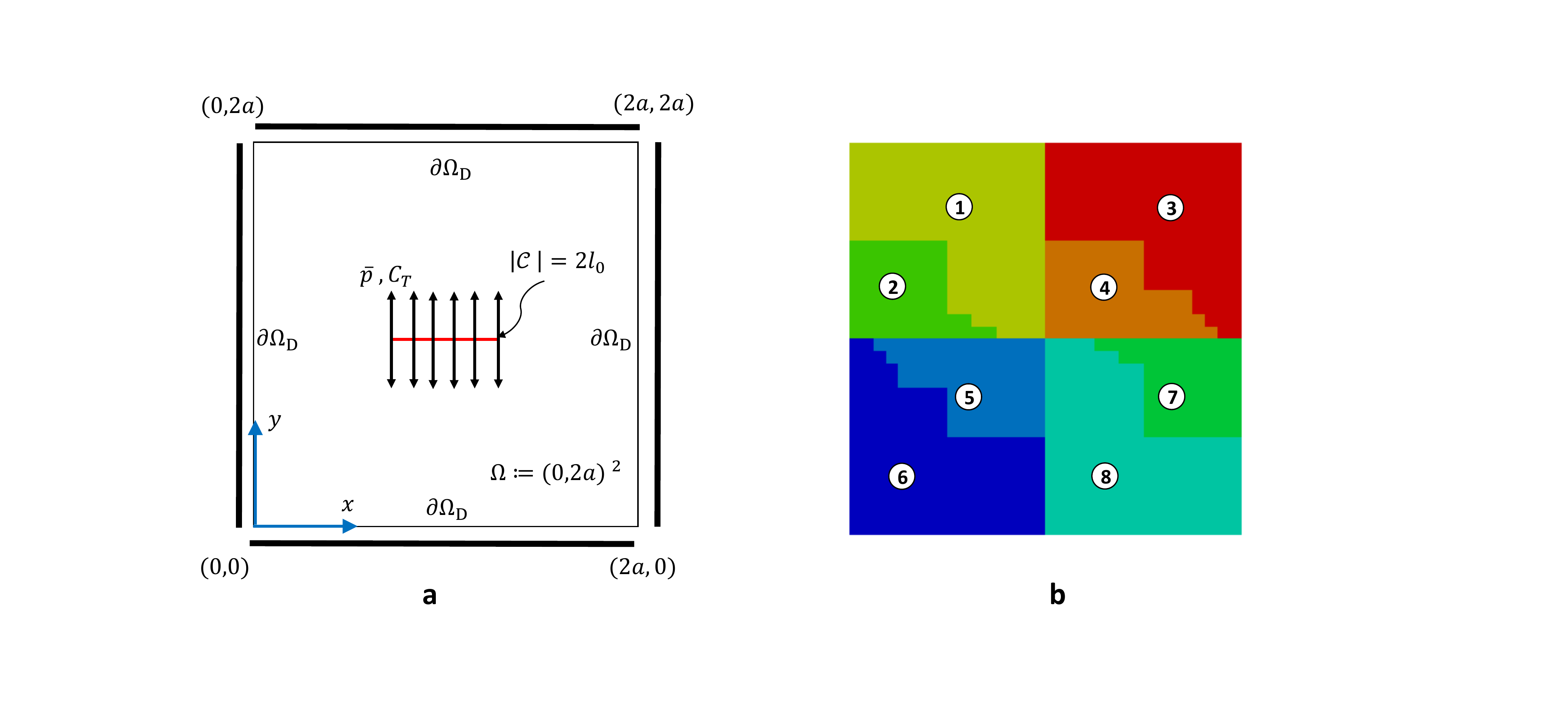}}   
	\caption{(a) Geometry and loading setup for the 2D test scenarios; (b) 
All numerical examples are computed by parallel computing on 8 processors. The different subdomains are associated with different processors. Depending on mesh refinement, the workload for each processor is adjusted dynamically at each time step.}
	\label{Figure3}
\end{figure}

In the first five examples, we present two-dimensional settings. The
geometrical parameter denoted as $a$ in Fig. \ref{Figure3} is set to 100, and,
hence a reservoir of size is $(0,200)^2\;m$. The Tran et al. problem  is
considered in an infinite domain with $\mathcal{C}$ with a pre-defined crack
of length $2l_0$ in the $y=a$ plane and is restricted in $a-l_0\le
|\mathcal{C}|\le a+ l_0$. We set the half crack length as a $l_0=10\;m$. 
We set the initial values for displacement and phase-field as $\bm u_0:=0 \in
\Omega$ and $\varphi_0:=1 \in \Omega_R$ and $ \varphi_0:=0 \in \Omega_F$.
The (finite) computational domain is subdivided into quadrilateral element
domains. 
It is noted, that all numerical examples are computed by parallel computing on
$8$ processors; see Fig. \ref{Figure3}[b]. The different sub-domains are associated with
different processors. Depending on mesh refinement, the workload for each
processor is adjusted dynamically at each time step. Typically, we set $2\times10^4$
quadrilateral elements for each processor.

\subsubsection{Three-dimensional setup}
\begin{figure}[H]
	\centering
	{\includegraphics[trim=0cm 0cm 0cm 0cm, width=13cm] {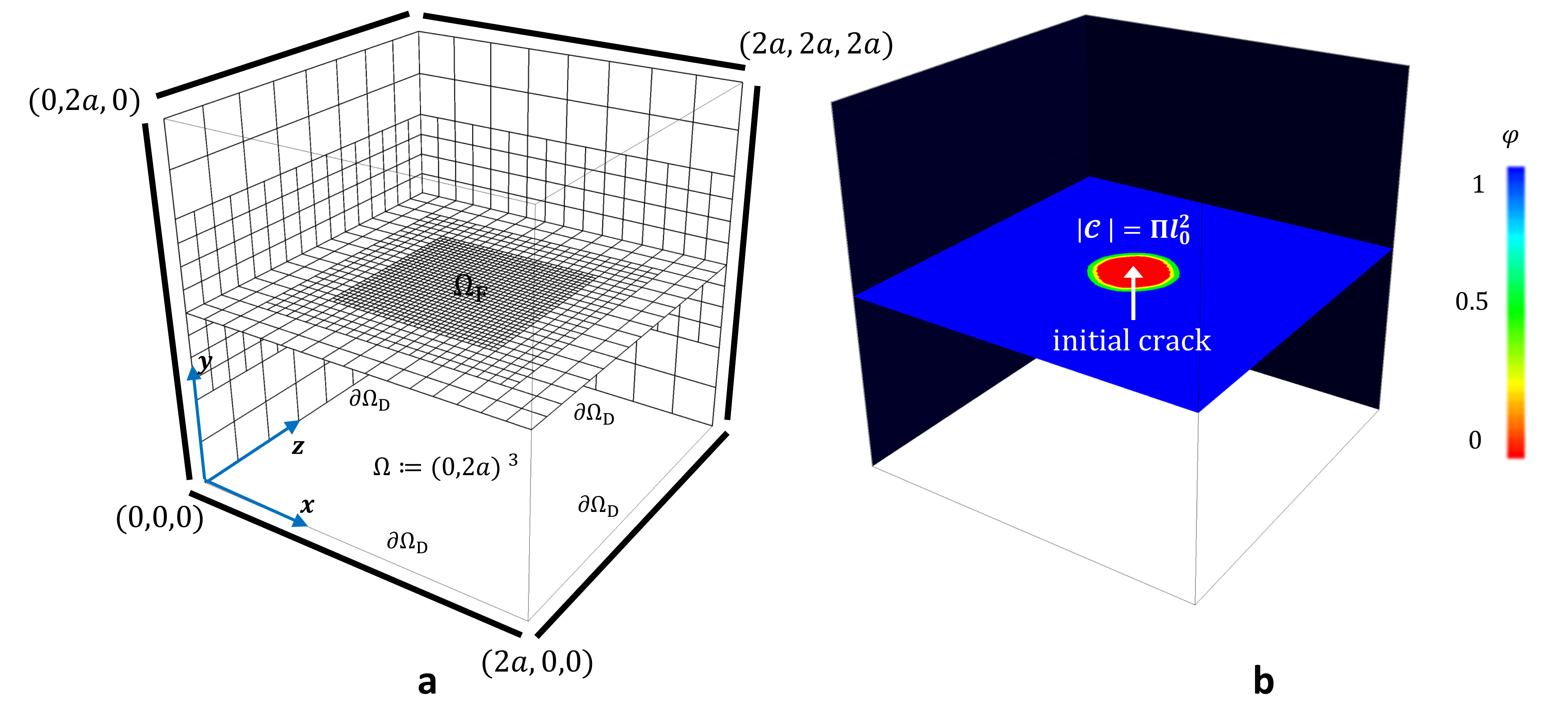}}
	\caption{Geometry and initial fracture for the 3D test scenarios. (a) locally refined mesh and (b) initial penny-shape crack in the middle.}
	\label{ex_2}
\end{figure}

In the last two examples, we present three-dimensional test cases. The
configuration setup is displayed in Fig. \ref{ex_2} and the values of the
material parameters are the same as in the two-dimensional test cases. The
geometrical parameter denoted as $a$ in Fig. \ref{ex_2} is set to 50, hence,
the reservoir of size is a $(0,100)^3\;m$ cube. We consider an infinite domain
with $\mathcal{C}$ as a pre-defined crack with a given radius $l_0$ in the
$y=a$ plane and is restricted in  $\mathcal{C}(x,y):= \{
(x,y)\in\Omega:(x-a)^2+(z-a)^2\leq l^2_0\}$. We set the crack radius length as
a $l_0=10\;m$. The (finite) computational domain is subdivided into hexahedral
element domains. 

\subsubsection{Material parameters}
\label{sec_mat_para}
The material parameters are taken from \cite{TrSeNg13} and related to
practical field problems.
The geo-mechanical parameters are Young's modulus $E_Y=1.5\times
10^{10}Pa$, Poisson's ratio $\nu = 0.15$, and the linear thermal-expansion
coefficient $\beta = 10^{-5} 1/C$ and finally the critical energy release rate
$G_c = 10^{10}\; N/m$.

For the flow and temperature, we have the thermal diffusivity
$\kappa_\Theta = 10^{-6} m^2/s$, the initial stress $p_0 = 12130\:$KPa and in
the stationary-in-length test cases (only evolution of the aperture), we use
the constant pressure  $\bar{p} = 15834\:$KPa. Additionally, for the
non-isothermal setting, the initial temperature $\Theta_0 = 100 \; C$, the
apparent temperature is set to $\Theta = 70 \: C$.

For all examples, we set $\alpha_B=0$ and $\alpha_\Theta=0$. 
All material properties are fixed for the following numerical examples, 
unless indicated otherwise.

\subsubsection{Model parameters}
	\label{Rem_heps}
The phase-field parameters are chosen as $\kappa = 10^{-10}$, 
and $\ep = 10\times \sqrt{h}$ (respecting the condition $h < \eps$) 
in 2D and $\ep = 2h$ in the 3D test cases.
	Using the adaptive mesh refinement, we are not dealing with a uniform
        mesh and hence the domain is divided into coarser and finer mesh
        elements, i.e. $\Omega=\Omega_c\cup \Omega_f$. Let
        $h^{max}_c:=max(h_c)$ in $\Omega_c$ and  $h^{min}_f:=min(h_f)$ in
        $\Omega_f$. On the one hand $\eps$,  enters in the constitutive
        modeling (that is in our PDE). On the other hand, it implicitly
        depends on the discretization of a domain, i.e. $h^{max}_c=o(\eps_c)$
        and $h^{min}_f=o(\eps_f)$. Thus we define, $\eps:=max(\eps_c,\eps_f)$
        resulting in $\eps\geqslant h_c>h_f$ and thus $\eps > h$ holds in
        every point of the domain. 
        
        The stopping criterion of the Newton method, i.e. the relative
        residual norm that is $\texttt{Residual}: \| \bm F(\bm x_{k+1}) \|
        \leq \texttt{Tol}_\texttt{N-R} \| \bm F(\bm x_{k}) \|$, 
        is set to $\texttt{Tol}_\texttt{N-R}=10^{-10}$.

%

The threshold value for the local predictor-corrector mesh refinement 
is $\texttt{TOL}_{\varphi} = 0.9$. That is, we refine when
\[
\varphi(x) < \texttt{TOL}_{\varphi} \quad\text{for } x\in \Omega.
\]

\subsection{Test scenarios}

We propose the following numerical tests:

\begin{itemize}
	\item \textbf{Case a}. Two-dimensional problem with constant pressure and without  thermal effects. 
	\item \textbf{Case b}. Two-dimensional  problem with constant pressure and fixed Hagoort’s decline constant value. That is time-independent problem and considered to validate proposed formulation for fixed thermal effect. 
	\item \textbf{Case c}. Two-dimensional  problem with constant pressure and evolving Hagoort’s decline constant through time. That is a time-dependent problem and used to validate proposed formulation for a temperature variation in time.
	\item \textbf{Case d}. Observing crack propagation in the two-dimensional problem, higher temperature
	difference is considered. That is time-dependent problem with a constant pressure and evolving Hagoort’s decline constant through the time.
	\item \textbf{Case e}. Observing crack propagation in the two-dimensional  problem, evolving in time for pressure and Hagoort’s decline constant are
          considered. This is a time-dependent problem. 
\item \textbf{Case f}. Three-dimensional problem with constant pressure and a fixed Hagoort’s decline constant value.  
\item \textbf{Case g}. Three-dimensional problem with increasing pressure and Hagoort’s decline constant are
          considered.
\end{itemize}

\subsection{Quantities of interest}

\begin{figure}
	\centering
	{\includegraphics[trim=0cm 3cm 0cm 3cm, width=15cm] {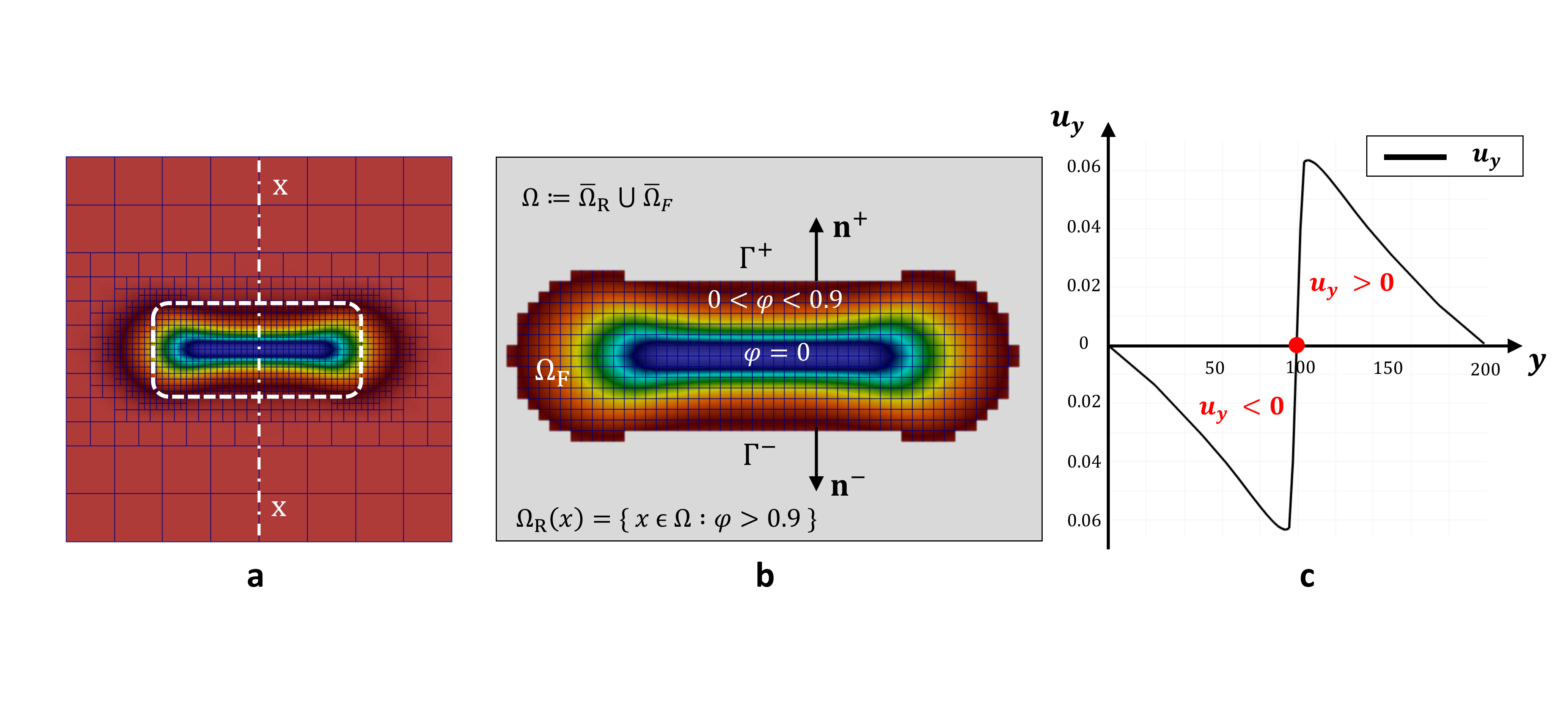}}   
	\caption{(a) Crack phase-field resolution (b) zoom into the framed region of the left plot and that is $\Omega_F:=\{x\in\Omega: \varphi(x,t)<0.9\}$ (c) Profile of the vertical deformation $\bm u_y$ on the section $x-x$ shown in (a) which shows symmetric displacements at the fracture boundary.}
	\label{Figure_disp_COD}
\end{figure}

In our numerical tests, we study the following aspects:
\begin{itemize}
\item Crack opening displacements (COD, also known as aperture) for the first
  three cases:
\begin{equation}
\label{COD}
COD := \int_{0}^{2a}[{\bm u}^+{\bm n}^+]-\int_{0}^{2a}[{\bm u}^-{\bm n}^-] \, dy=\int_{0}^{2a}[{\bm u}^++{\bm u}^-].\bm n \, dy=\int_{0}^{2a} \bm u(x_0, y)\cdot \nabla\varphi(x_0,y) \, dy,
\end{equation}
where ${\bm n}:={\bm n}^+=-{\bm n}^-$ is the normal vector (see
Fig. \ref{Figure_disp_COD} b) and $x_0$ the $x$-coordinate of the integration line.
We note that the integration is perpendicular to the crack direction.  Here, the crack is aligned with the $x$-axis and therefore integration into the normal direction coincides with the $y$-direction.
The analytical solution for the COD derived by Tran et al. \cite{TrSeNg13} reads:
\begin{equation}
\label{COD_analytical}
w(x,t) = \frac{2(1-\nu_s^2)l_0}{E} \sqrt{1-\rho^2}\;
\bigl( p-p_0 - C_\Theta(\Theta - \Theta_0) \bigr),
\end{equation}
where $l_0$ is the half-length of the crack and $0<\rho<1$, for $\rho = x/l_0$
with $x$ being the distance to the origin of the crack. In three dimensions,
the formula reads:
\begin{equation}
\label{COD_analytical_3D}
w (x,t) = \frac{4(1-\nu_s^2)l_0}{\pi E} \sqrt{1-\rho^2}
\bigl( p-p_0 - C_\Theta(\Theta - \Theta_0) \bigr).
\end{equation}
For pressurized fracture without any thermal effect, i.e., $\Theta = \Theta_0
= 0$, we obtain the formulas derived in Sneddon and Lowengrub
\cite{SneddLow69}[Section 2.4 and Section 3.3].

\item Fracture length and path for propagating fractures.
\item Both stopping criteria of the nonlinear solver (Section
  \ref{sec_primal_dual_as}), i.e. the relative residual norm and the active set constraint.
\item For some fixed time steps, the average number of GMRES iterations within
  per Newton cycle.
\item The evolution of the mechanical strain energy functional,
  i.e. mechanical term in Formulation \ref{form_3} and dissipated fracture
  energy functional, i.e. the fracture term in Formulation \ref{form_3}.
\item Comparisons regarding no split and vol./dev. split of strain density function.

\end{itemize}

%

\subsection{Case a. Sneddon-Lowengrub's 2D setting with pressure}
In the first example, we consider a stationary setting with a pressurized
fracture. The geometrical setup is the same as illustrated in Fig. \ref{Figure3} and 
the material parameters are listed in Section \ref{sec_mat_para}. 
Related results were
presented in \cite{MiWheWi18, WheWiWo14, BourChuYo12}.

We set the time step size $\delta t=1\:s $ for a given 
constant pressure $\bar{p}=10^{-3} \:KPa$. This test case is 
computed in a quasi-stationary manner: that is, we solve several 
pseudo-time steps. The goals of this test are to observe the crack 
opening displacement, and in particular to show the crack 
tip approximation. 
The initial mesh is five times uniformly refined, 
then, local mesh adaptivity is applied for five further levels. 

\begin{figure}[!ht]
	\centering
	{\includegraphics[clip,trim=0cm 2cm 0cm 5cm, width=15cm] {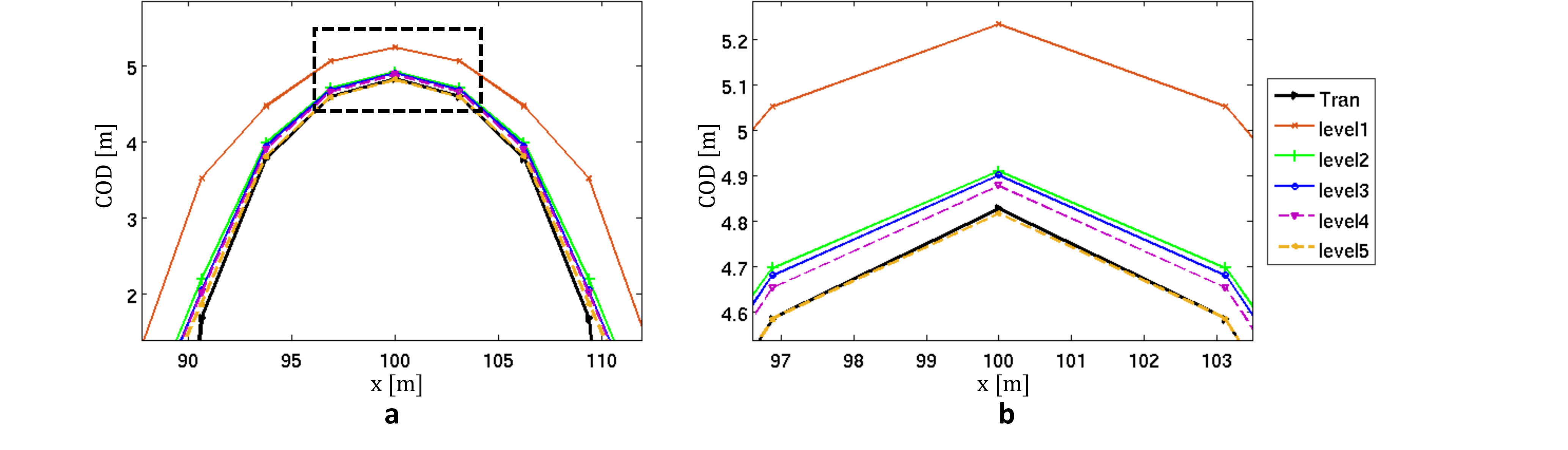}}   
	\caption{Case a. (a) COD for different global refinement level. (b)
          Zoom into the framed region on the left figure. Here the Tran et
          al. solution corresponds to the Sneddon-Lowengrub formula since no 
          temperature is considered in this first example.}
	\label{Figure4}
\end{figure}

The COD findings for both our computational model and the analytical
calculation are shown in Fig. \ref{Figure4}. It is observed the crack tips at
$90m$ and $110m$ show excellent agreement with the analytical solution and convergence towards Sneddon-Lowengrub’s manufactured  solution.  

\newpage
\subsection{Case b. 2D setting with pressure and constant temperature}
\begin{figure}[!ht]
	\centering
	{\includegraphics[clip,trim=0cm 2cm 0cm 4cm, width=15cm] {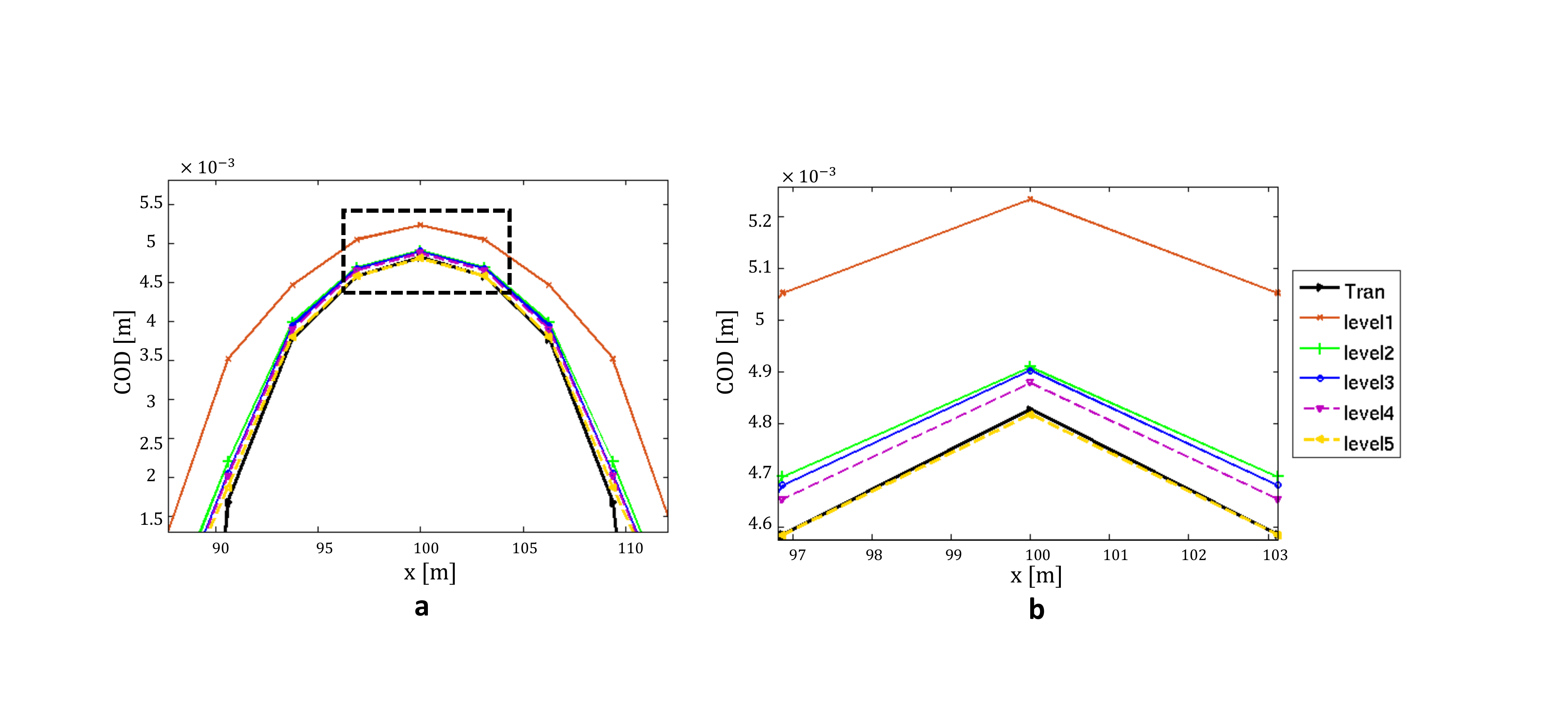}}   
	\caption{Case b. (a) COD for different refinement levels. (b) Zoom into the framed region on the left figure.}
	\label{Figure5}
\end{figure}

In this second example, we turn our attention now to a pressurized,
non-isothermal configuration. Specifically, we follow
Tran's et al. \cite{TrSeNg13} 2D problem with a
constant pressure and a fixed Hagoort's decline constant. We keep the geometry,
all parameters, and tolerances as in the first example. Additionally, we set
$\lambda_\Theta = 10^{-4}$ constant and we work still in a full stationary
setting where the time step loop terminates after one step. The crack opening
displacement both obtained through phase-field fracturing modeling and the 
analytical calculation are shown in Fig. \ref{Figure5} for 
the same refinement levels as in the previous example.

\begin{figure}[!ht]
	\centering
	{\includegraphics[clip,trim=0cm 0cm 0cm 0cm, width=15cm] {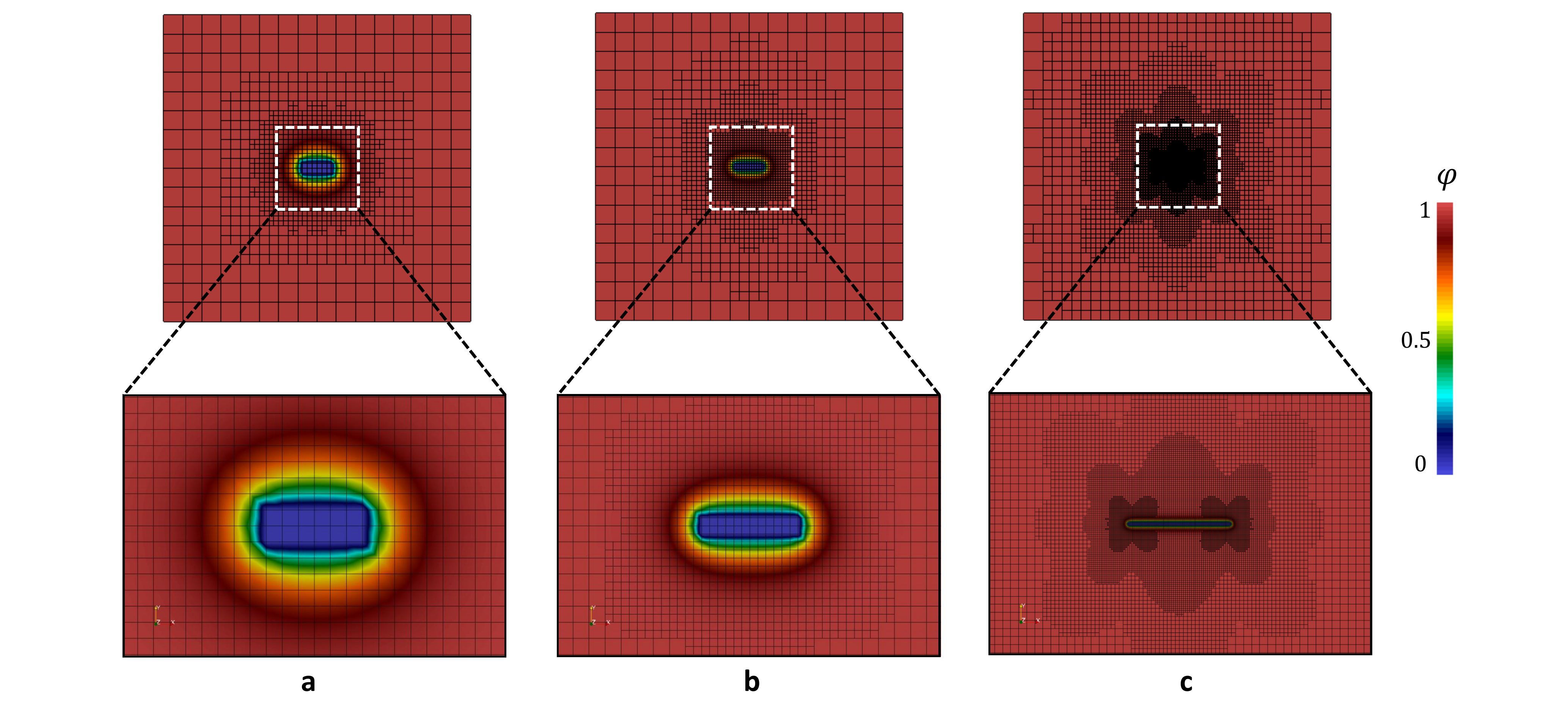}}   
	\caption{Case b. Crack phase-field resolution in the range of
          different time steps at time (a) 9 sec. (b) 14 sec. (c) 25 sec. on
          different locally refined meshes. Here, the fracture only varies 
in its width (i.e., COD), but not in length. Nonetheless, we perform time step
solutions
in order to fulfill the crack irreversibility constraint.}
	\label{Figure6}
\end{figure}

The crack phase-field pattern for different time steps on the locally adaptive
refined meshes are illustrated in Fig. \ref{Figure6}.

\subsection{Case c. 2D setting with pressure, temperature and decline constant}
\label{sec_case_c}
While Hagoort's decline constants were constant in Case b, we now 
account for a temperature variation 
in time. We keep the geometry, all parameters, 
and tolerances as in the first example. The only 
two changes are that we set 
the time step size to one day $\delta t=86400s$ and 
observe $T = 365$ days. 
This test case has a slight time-dependence since
we use $\lambda_\Theta :=\lambda_\Theta(t)$ being time-dependent
and defined through Eq.  \ref{L_T_eq}. 
This scenario allows us to study the fracture behavior 
for a constant pressure like in the production processes.

\begin{figure}[!ht]
	\centering
	{\includegraphics[clip,trim=0cm 3cm 0cm 3cm, width=15cm] {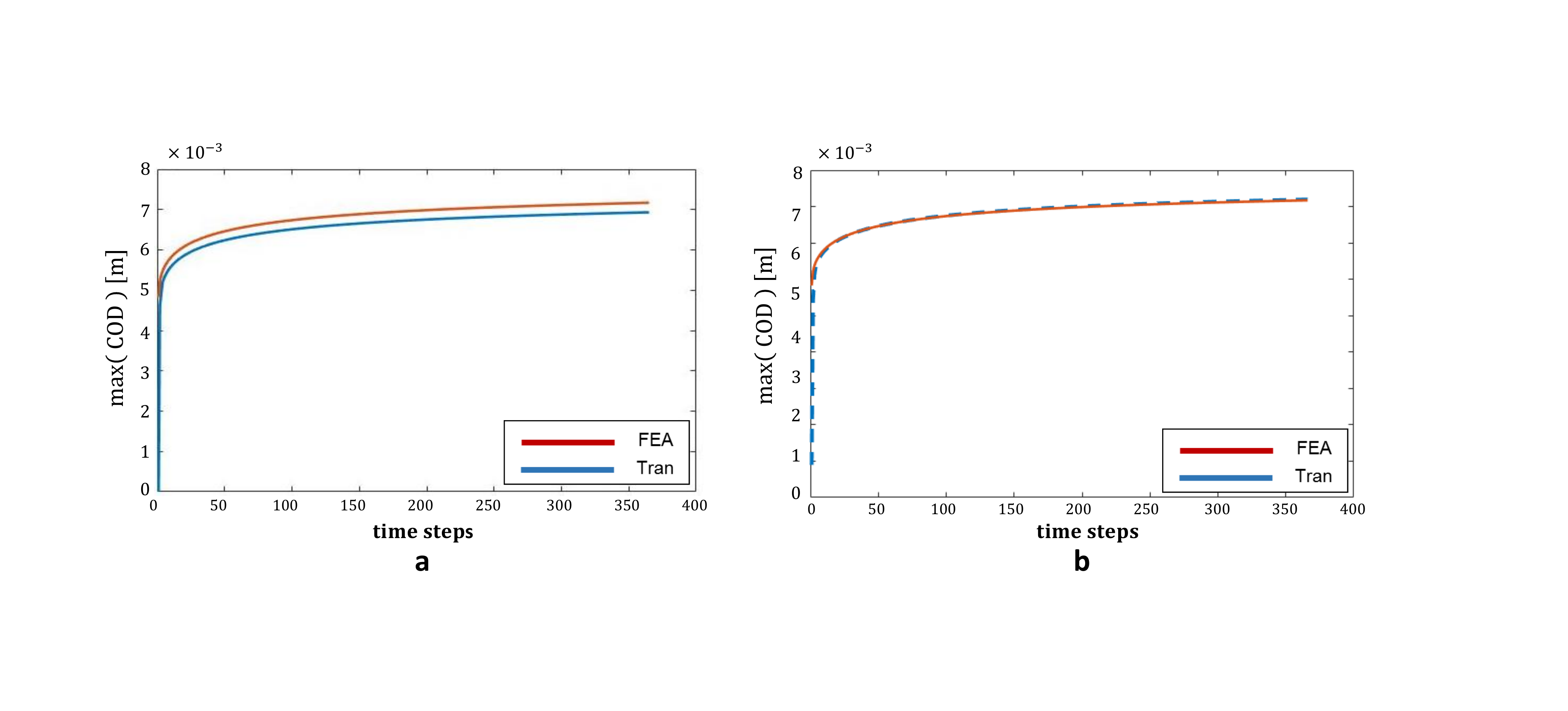}}   
	\caption{Case c. Maximum width evolution over time in the 
		non-isothermal case. Injection of cold water $\Theta < \Theta_0$ leads 
		to an increasing aperture over large time scales. (a) level 5 and (b) level 9 global refinement scheme.}
	\label{Figure7}
\end{figure}

Our numerical findings regarding the maximum COD is shown 
in Fig. \ref{Figure7}. 
We consider different levels for uniformly 
refined meshes. Findings for level 5 and level 9 are shown in Fig. \ref{Figure7} a and b, 
respectively. In between we computed the levels 6-8 as well, but are not shown 
to keep the figures legibly. 
We observe excellent agreement of our numerical model with 
the analytical solution.

\begin{figure}[!ht]
	\centering
	{\includegraphics[clip,trim=0cm 0cm 0cm 0cm, width=15cm] {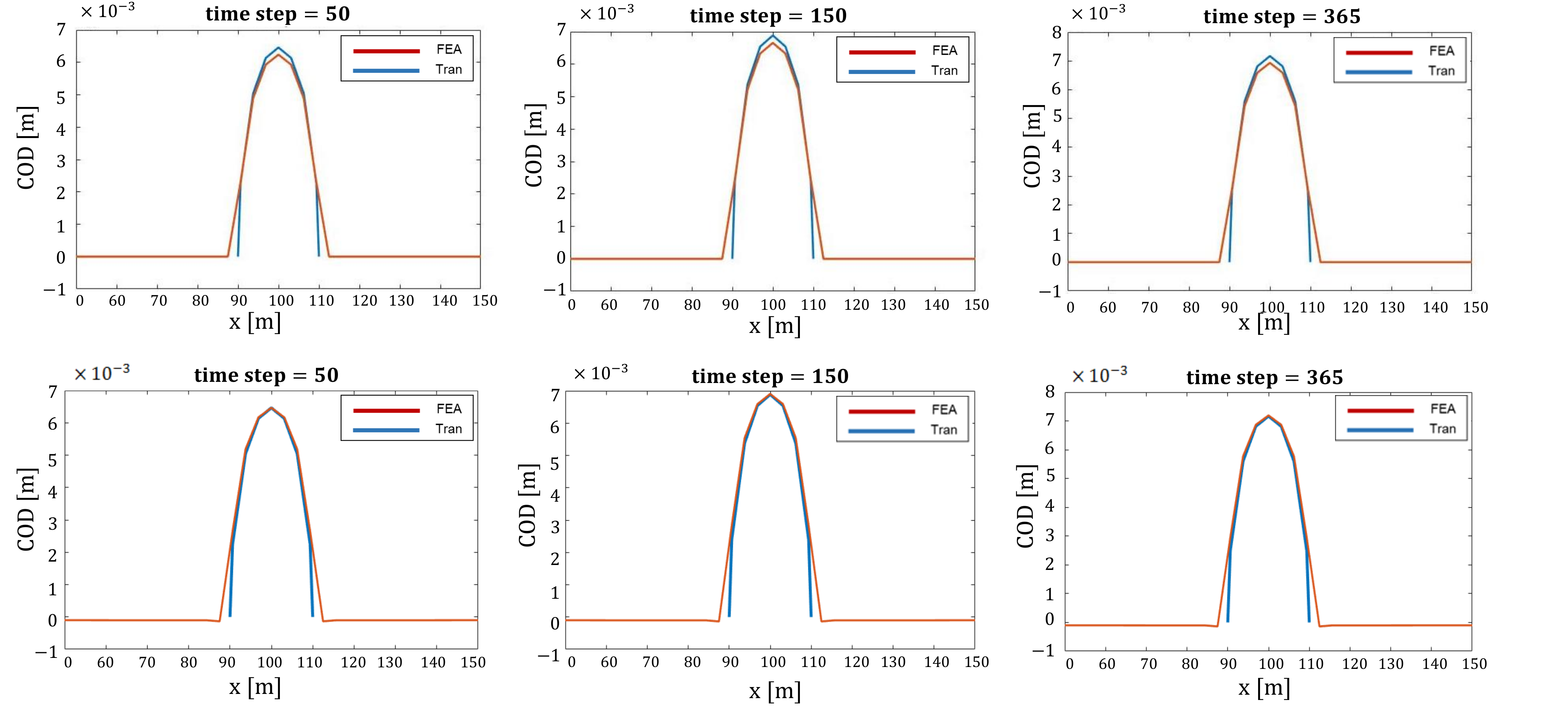}}   
	\caption{Case c. Comparison of the COD between phase-field numerical method and Tran's et al. \cite{TrSeNg13} analytical solution during different time steps (a) first row corresponds to the level 5 and (b) second row corresponds to the level 9 global refinement scheme.}
	\label{Figure8}
\end{figure}

\newpage
Another perception is to observe again the COD over the fracture (similar to
Case a and Case b). This can be shown by fixing the time steps 
(i.e. 50, 150 and 365 days). Our findings are illustrated in
Fig. \ref{Figure8} and show again excellent agreement.

\subsection{Case d. 2D setting for the crack growth due to temperature variation}
\begin{figure}[!ht]
	\centering
	{\includegraphics[clip,trim=0cm 2cm 0cm 1cm, width=15cm]
          {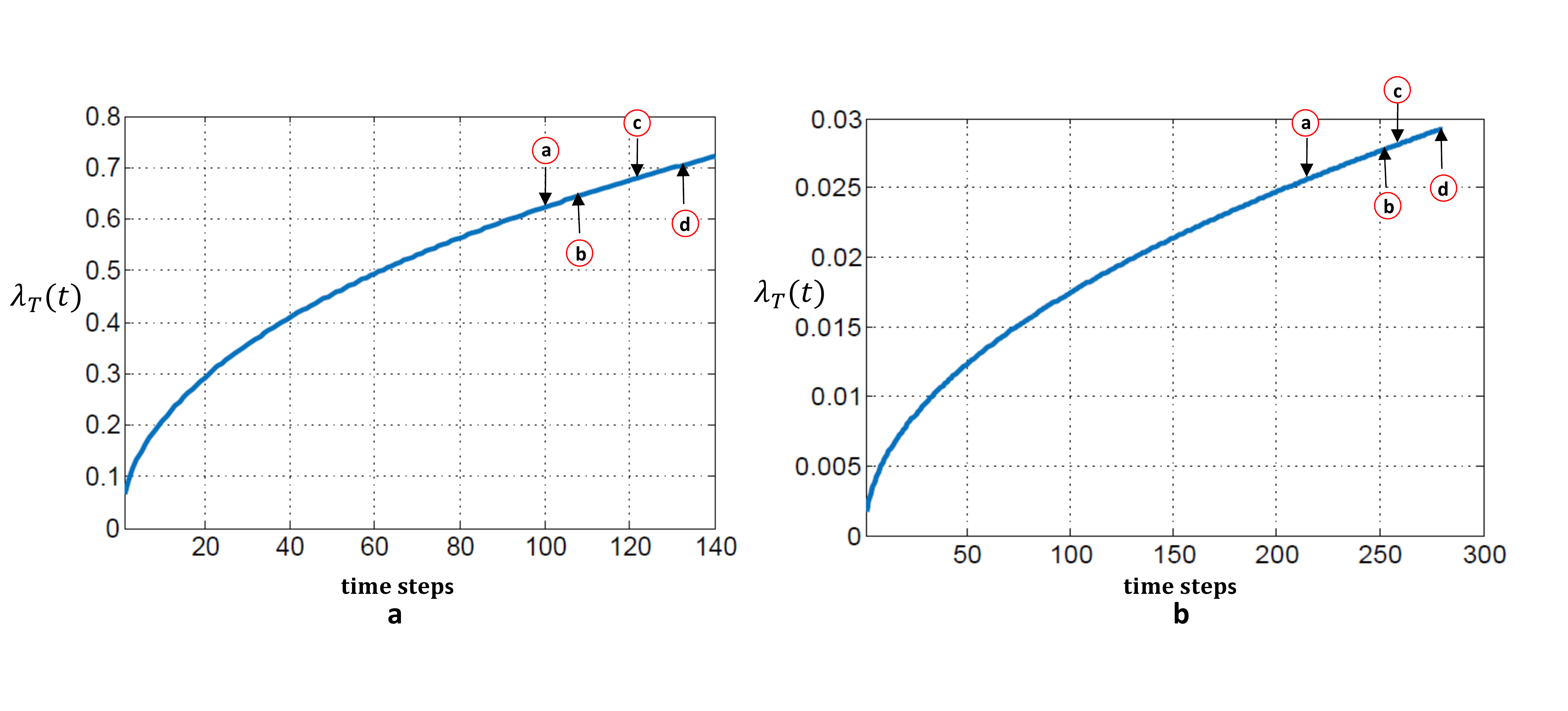}}   
\caption{Hagoort's decline constants for a time-dependent problem that is $\lambda_\Theta :=\lambda_\Theta(t)$, i.e. Eq. \ref{L_T_eq}, for (a) Case d and (b) Case e. Interested time step steps for each cases are shown in each plot.}
	\label{Figure9}
\end{figure}

Having studied slight time dependencies with variations 
in the COD, we now consider a first test case with a propagating fracture.
The propagation will be caused by temperature variations.
For temperature variations, we use
$\lambda_\Theta :=\lambda_\Theta(t)$. Hagoort's decline constants for this
example are shown in Fig. \ref{Figure9} a. 
The second novelty in this numerical test are comparisons of 
strain-energy split and the non-split version.
The third goal are very detailed studies on the performance 
of the linear and nonlinear solvers.

Again, we keep most
parameters as in the first example except a higher temperature difference to $220$
degree Celsius (i.e. $\Theta_0=300\: C \; \text{and} \; \Theta=80\: C$) and the
initial stress $p_0 = 12130\:$KPa and constant pressure as mentioned earlier
is set to $\bar{p} = 15834\:$KPa. 

In order to facilitate fracture propagation, 
the critical energy release rate is
drastically reduced to $G_c = 5.5\times 10^5\; N/m$.
The time step size is $\delta t=86400s$. 
As it is shown in Fig. \ref{Figure9} a, we are now interested 
in four specific time steps, namely $100, 114, 121$ and $136$ days.

\subsubsection{Crack propagation and locally refined meshes}

\begin{figure}[!ht]
	\centering
	{\includegraphics[clip,trim=0cm 0cm 0cm 0cm, width=15cm] {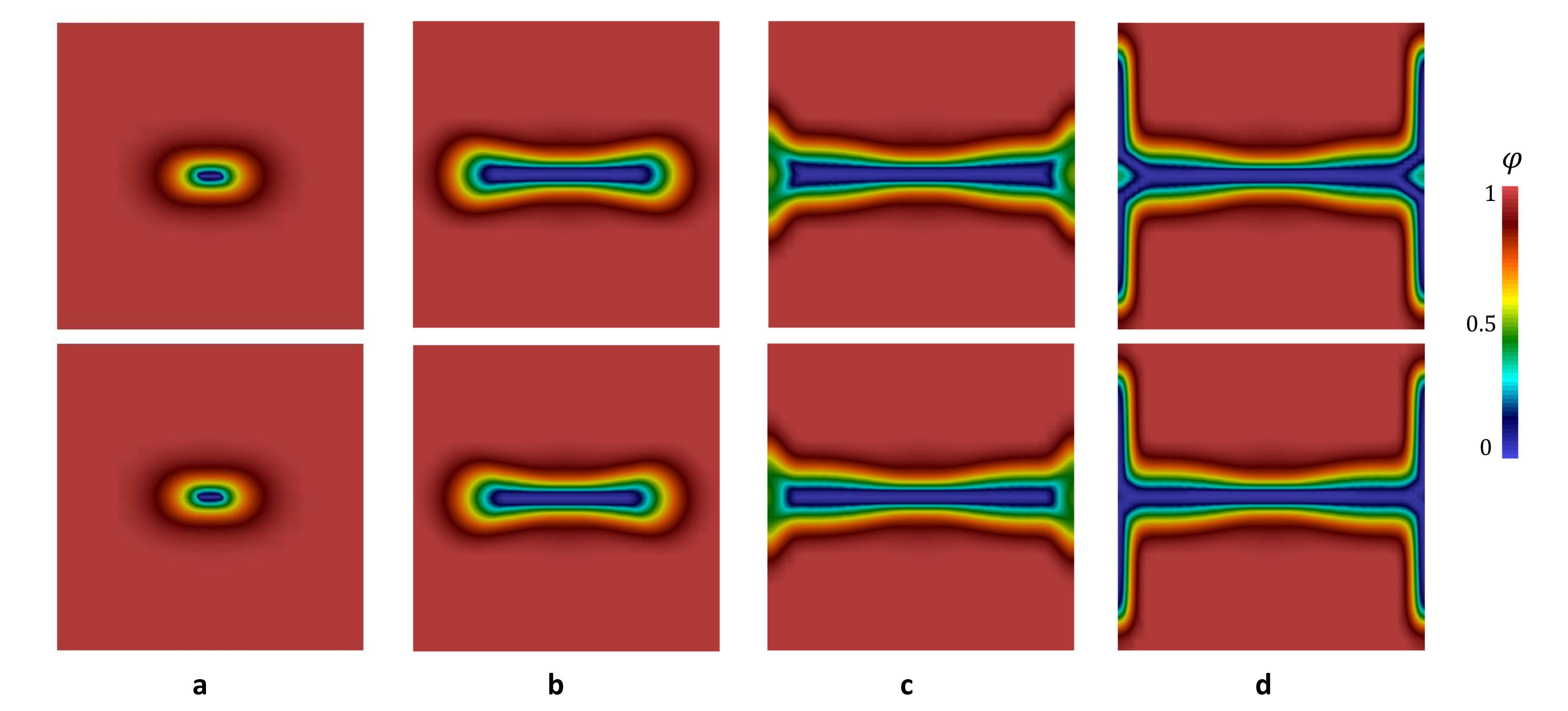}}   
	\caption{Case d. Crack phase-field evolution due to temperature variations at (a) 100 day (i.e $\lambda_\Theta= 0.6225 \: ,C_\Theta=48933$) (b) 114 day ($\lambda_\Theta= 0.6596 \: ,C_\Theta=50189$) (c) 121 day ($\lambda_\Theta= 0.6770 \: ,C_\Theta=50751 $) (d)136 day (i.e $\lambda_\Theta= 0.7121 \: ,C_\Theta=51837$). First and second row represented for the no split and vol./dev. split of the strain energy density function, respectively.}
	\label{Figure10}
\end{figure}

Figure \ref{Figure10} displays a sequence of crack 
patterns during different time steps without considering 
any strain-energy split, i.e. Formulation \ref{form_2} and considering vol./dev. split, 
i.e. Formulation \ref{form_3}, of the strain energy density function. 
In fact, until time step $100$, we observe almost no growth. Then, 
we have crack propagation towards the boundaries.  

\begin{figure}[!ht]
	\centering
	{\includegraphics[clip,trim=1cm 5cm 0cm 3cm, width=15cm] {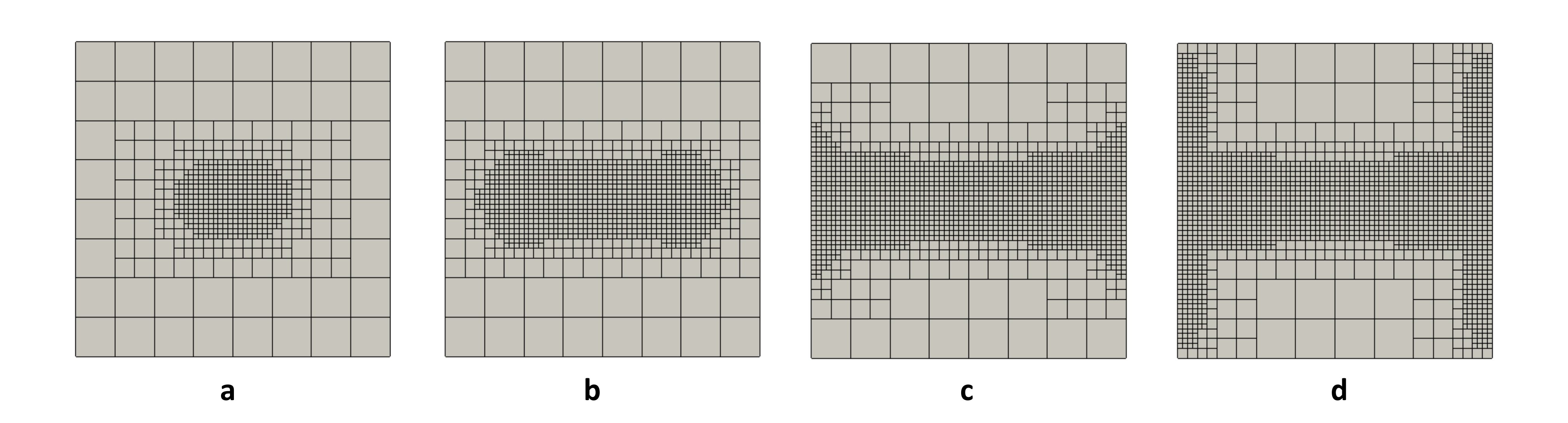}}   
	\caption{Case d. Functionality of predictor-corrector mesh refinement, i.e. the mesh evolves
		with the fracture. The transition zone with $0<\varphi<1$ determines the region in which the
		mesh has to be refined. $\texttt{Tol}_{\varphi}$ is set to 0.9. The refine meshes represent times steps (a) 100 day (b) 114 day (c) 121 day (d)136 day.}
	\label{Figure11}
\end{figure}

\newpage
The functionality of predictor-corrector mesh refinement 
is shown in  Fig. \ref{Figure11} to present the evolution of 
the locally refined mesh when the crack is growing. 
As mentioned earlier in 
Section \ref{Rem_heps}, $\mathrm{Tol}_{\varphi}$ is set to 0.9.
Each cell that has 
at least one support point with value $\varphi(x_i)<\texttt{Tol}_{\varphi}$ 
for $x_i\in \Omega_e$ will be refined unless we are already at the maximum 
desired refinement level is reached.

\subsubsection{Analysis of the strain-energy splitting}
We now turn to the second goal that is novel in this form 
in the published literature regarding pressurized/non-isothermal fracture
settings with phase-field modeling.

\begin{figure}[!ht]
	\centering
	{\includegraphics[clip,trim=0cm 0cm 0cm 0cm, width=15cm] {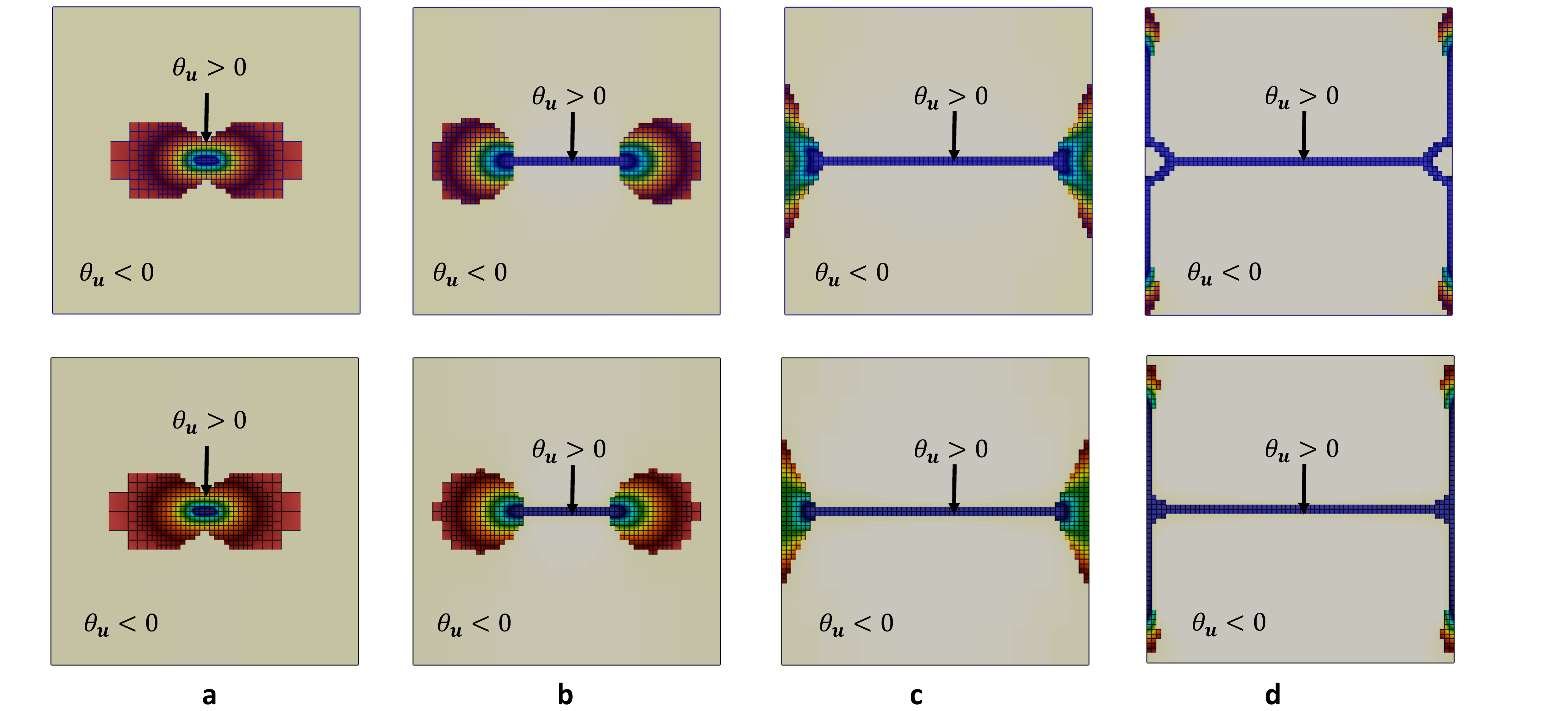}}   
	\caption{Case d. Qualitative representation of the tensile and
          compression counterpart of domain. The compression region,
          i.e. $\theta_{\bm u}:=H{^+}(\nabla. \bm u)<0$, is shown uniformly in
          gray, and the tensile region, i.e. $\theta_{\bm u}>0$, is shown with
          desaturated rainbow colors at (a) 100 day (b) 114 day (c) 121 day
          (d) 136 day. The first row represents no splitting and the second 
          row shows split of the strain energy density function.}
	\label{Figure12}
\end{figure}

A qualitative representation of the tensile and compression counterparts 
are shown in Fig. \ref{Figure12}. We define 
$\theta_{\bm u}:=H{^+}(\nabla. \bm u)=H{^+}(tr(\bm \varepsilon))$ that is an indicator 
to detect the compression region, i.e. $\theta_{\bm u}<0$, 
and, reversily, to detect the tensile region, i.e. $\theta_{\bm u}>0$. 
If the vol./dev. split of the strain energy density function is used, 
we are in the crack region when the phase-field energy exceeds its critical 
value $G_c$ and additionally $\theta_{\bm u}>0$. This gives us an additional
constraint to our physical model. Hence, only the tensile stress is degraded
and compression stress is 
free from the effect of the phase-field variable; see Eq. \ref{Q7}.

\begin{figure}[!ht]
	\centering
	{\includegraphics[clip,trim=0cm 1cm 0cm 3cm, width=15cm] {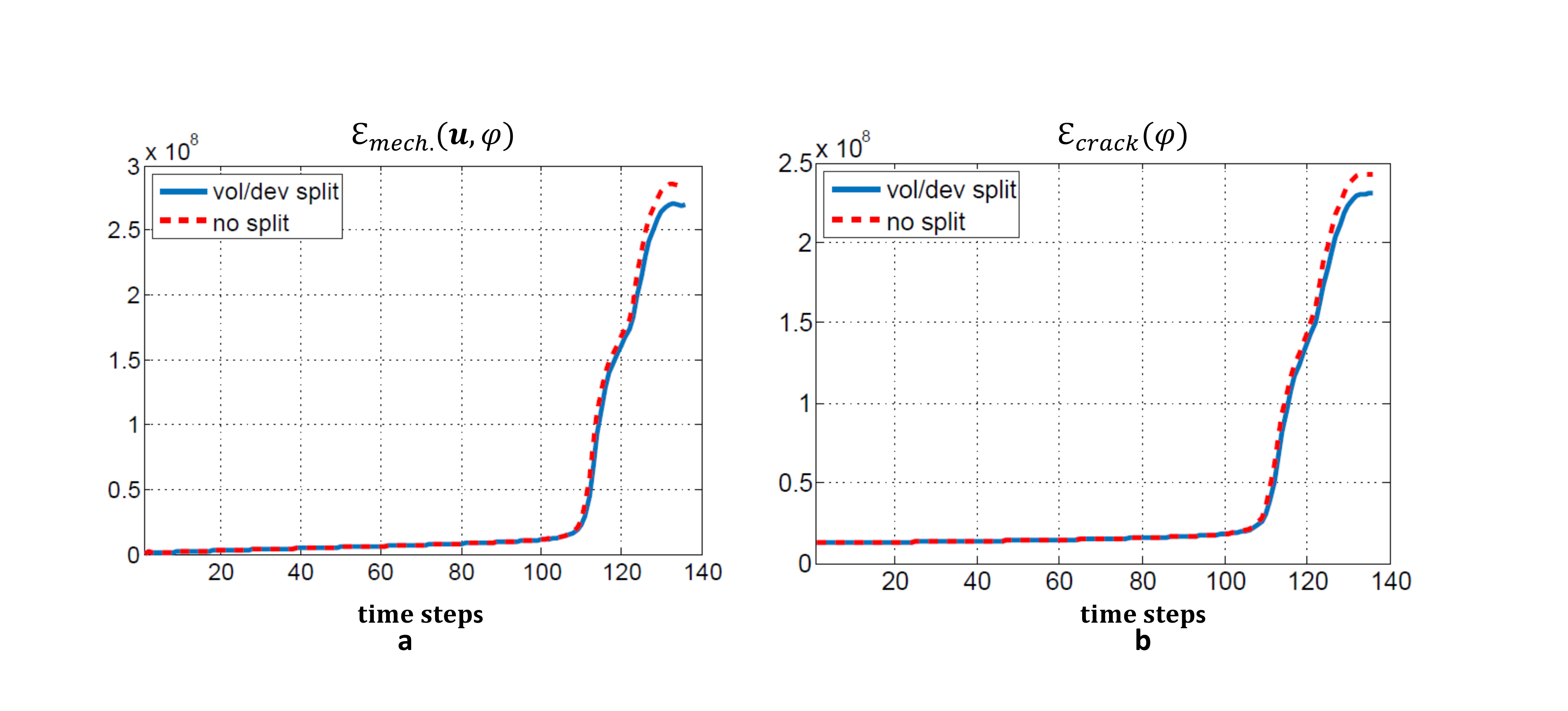}}   
	\caption{Case d. Plots over time of (a) the elastic strain energy per unit thickness, i.e. mechanical term in formulation \ref{form_3} (b) the dissipated fracture energy per unit thickness , i.e. fracture term in formulation \ref{form_3}, with no split and vol./dev. split for strain energy density function.}
	\label{Figure13}
\end{figure}

Figure \ref{Figure13} compares the incremental elastic strain energy, i.e. 
the mechanical term in the Formulation \ref{form_3}, and the crack dissipated energy, 
i.e. fracture term in Formulation \ref{form_3}. 
For comparison, the strain energy density function without 
split and vol./dev. split are shown. As observed in Fig. \ref{Figure13}, it turns 
out that dissipated fracture energy resulting from the crack phase-field 
evolution evolve with a sharp ascending behavior.  This clearly shows that the 
sharp transition between the intact state and the fully cracked state of the solid 
body arised. Thereafter, dissipated fracture energy evolves towards a steady-state 
response.  While the solid body completely fractured, no change in 
the dissipated fracture energy is observed anymore. 

\subsubsection{Performance of nonlinear and linear solvers}

\begin{figure}[!ht]
	\centering
	{\includegraphics[clip,trim=0cm 0cm 0cm 0cm, width=15cm] {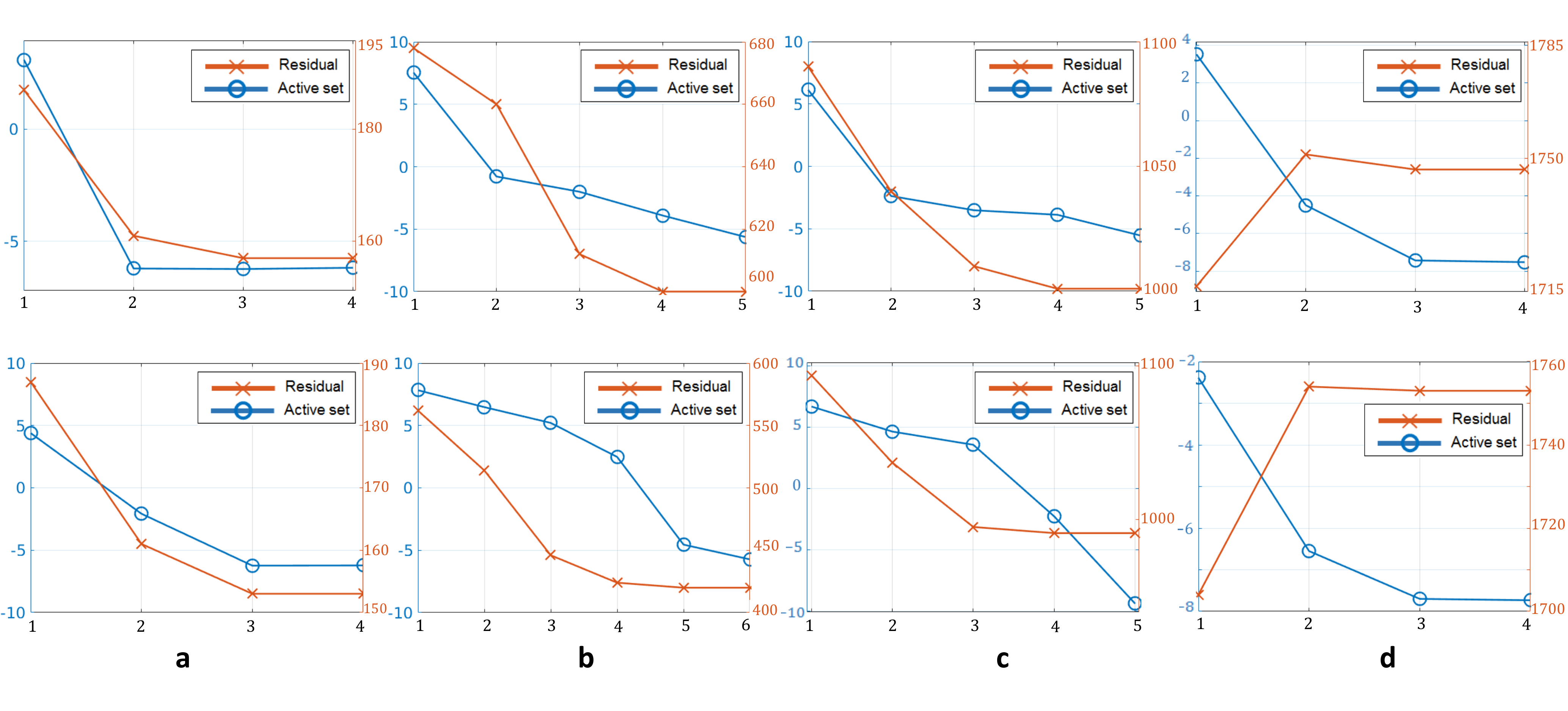}}   
	\caption{Case d. Performance of the combined Newton
          solver. Convergence behavior for the nonlinear relative residual,
          i.e. red line, and elements in the active set, i.e. blue line, at
          the fixed time steps (a) 100 day (b) 114 day (c) 121 day (d)136
          day. The first row represents no splitting and the second 
          row shows split of the strain energy density function.
}
	\label{Figure131}
\end{figure}

A detail analysis of the convergence of the primal-dual active set
framework and  the Newton solver are shown in Fig. \ref{Figure131} for 
the finest mesh level. 

Finally, the average number of GMRES iterations of linear iterations per Newton cycle for 
different times steps is shown in Table \ref{GMRES_itr_caseD}. 
Specifically, the preconditioner (Section \ref{sec_lin_sol}) works extremely well.

\begin{table}[!ht]\small
	\caption{Case d. Average number of GMRES iterations for one Newton
          cycle for four different times steps.}
	\label{GMRES_itr_caseD}
{
	\begin{center}
	\begin{tabular}{c|c|c}
		\hline
		Time step (day) & no split &vol./ dev. split\\
		\hline 
		100 & 2 & 2 \\ 
		114 & 6 & 15\\
		121 & 6 & 16\\
		136 & 16 & 17\\ 
		\hline
	\end{tabular}
	\end{center}
}
\end{table}

\subsection{Case e. 2D setting for pressurized crack propagation in non-isothermal settings}
\label{sec_case_e}
In this example, we extend the previous Case d to study now crack propagation
due to temperature effects , i.e. $\lambda_\Theta :=\lambda_\Theta(t)$ and
increasing pressure, i.e. $p:=p(t)$.  
From the application viewpoint this example is closer to a fracture process 
with a characteristic time scale much smaller than in Case d. 
Hence, we set the time step size to one day 
$\delta t=60s$ and simulate $270$ minutes.  The increasing pressure is prescribed as:
\begin{equation}
p:=p(t) = \alpha(t+1)\times \bar{p}=\alpha(t+1)1.5834\times 10^{7}, \quad \text{for} \quad t\geq0.
\end{equation}
We set the constant value $\alpha:=\frac{1}{2}$ and $t$ is the 
time and given by $t=n\delta t$ with $n$ being the time 
step number. We keep all remaining parameters as in the first example and we set $\Theta_0=100\; C$ 
and $\Theta=70 \; C$ degree Celsius (recall $\Theta< \Theta_0$, see Remark \ref{Rem1}).

\subsubsection{Crack propagation and locally refined meshes}

\begin{figure}[!ht]
	\centering
	{\includegraphics[clip,trim=0cm 0cm 0cm 0cm, width=15cm] {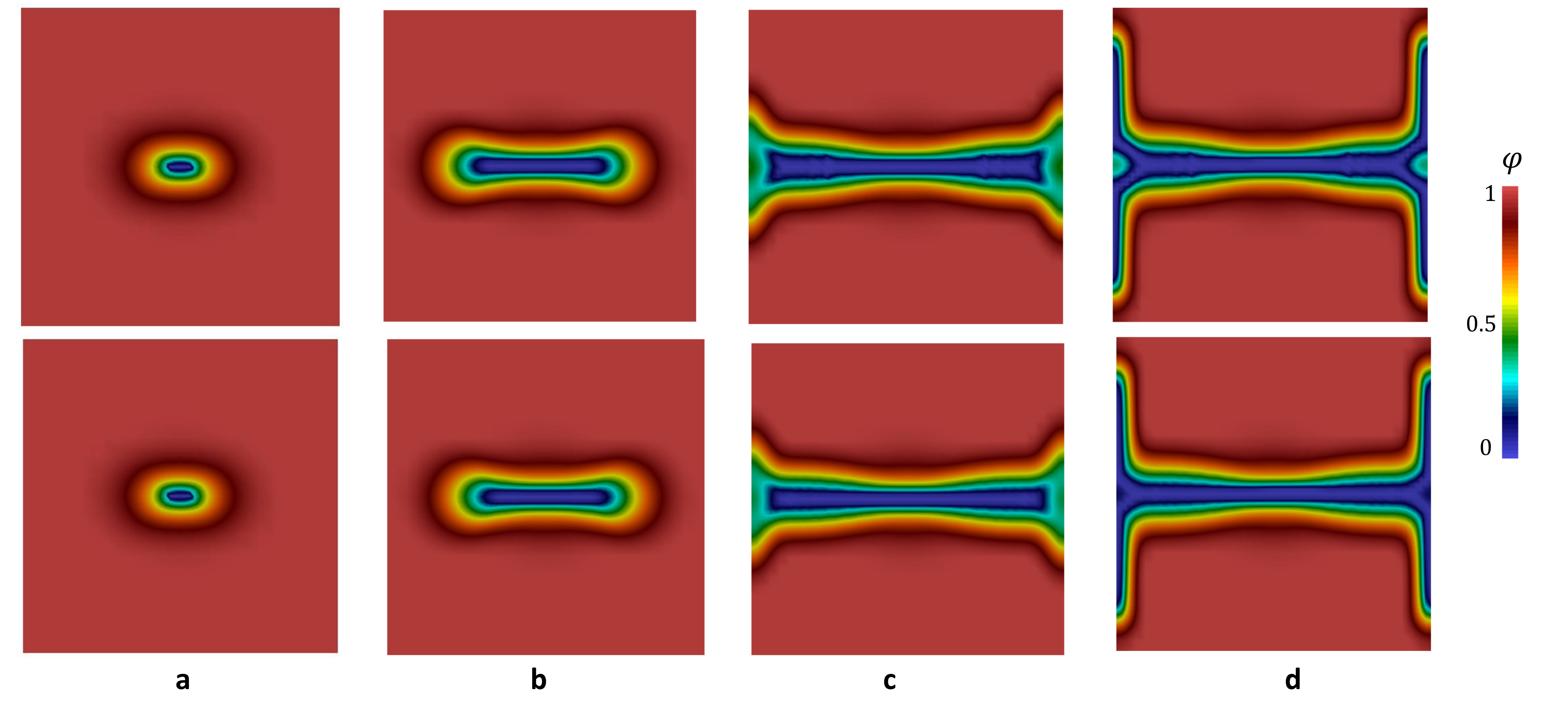}}   
	\caption{Case e. Crack phase-field evolution due to temperature and
          pressure variations at (a) 240 min (i.e. $\lambda_\Theta= 0.0271 ,
          C_\Theta= 4534, p=1.9080\times 10^9$) (b) 252 min
          (i.e. $\lambda_\Theta= 0.0278, C_\Theta= 4640.1, p=2.0030\times
          10^9$) (c) 259 min (i.e. $\lambda_\Theta= 0.0281, C_\Theta= 4700.7,
          p=2.0584\times 10^9$)  (d) 267 min (i.e. $\lambda_\Theta= 0.0286,
          C_\Theta= 4768.8, p=2.1218\times 10^9$). 
The first row represents no splitting and the second 
          row shows split of the strain energy density function.}
	\label{Figure14}
\end{figure}

As shown in Fig. \ref{Figure14} b, 
we are now interested in four time steps, namely 240, 252, 259 and 267 minutes
to observe the crack path. There again (as in Case d) only slight differences 
at the same points, suggesting that strain-splitting has nearly no effect 
in these configurations.

\subsubsection{Analysis of the strain-energy splitting}

\begin{figure}[!ht]
	\centering
	{\includegraphics[clip,trim=0cm 0cm 0cm 0cm, width=15cm] {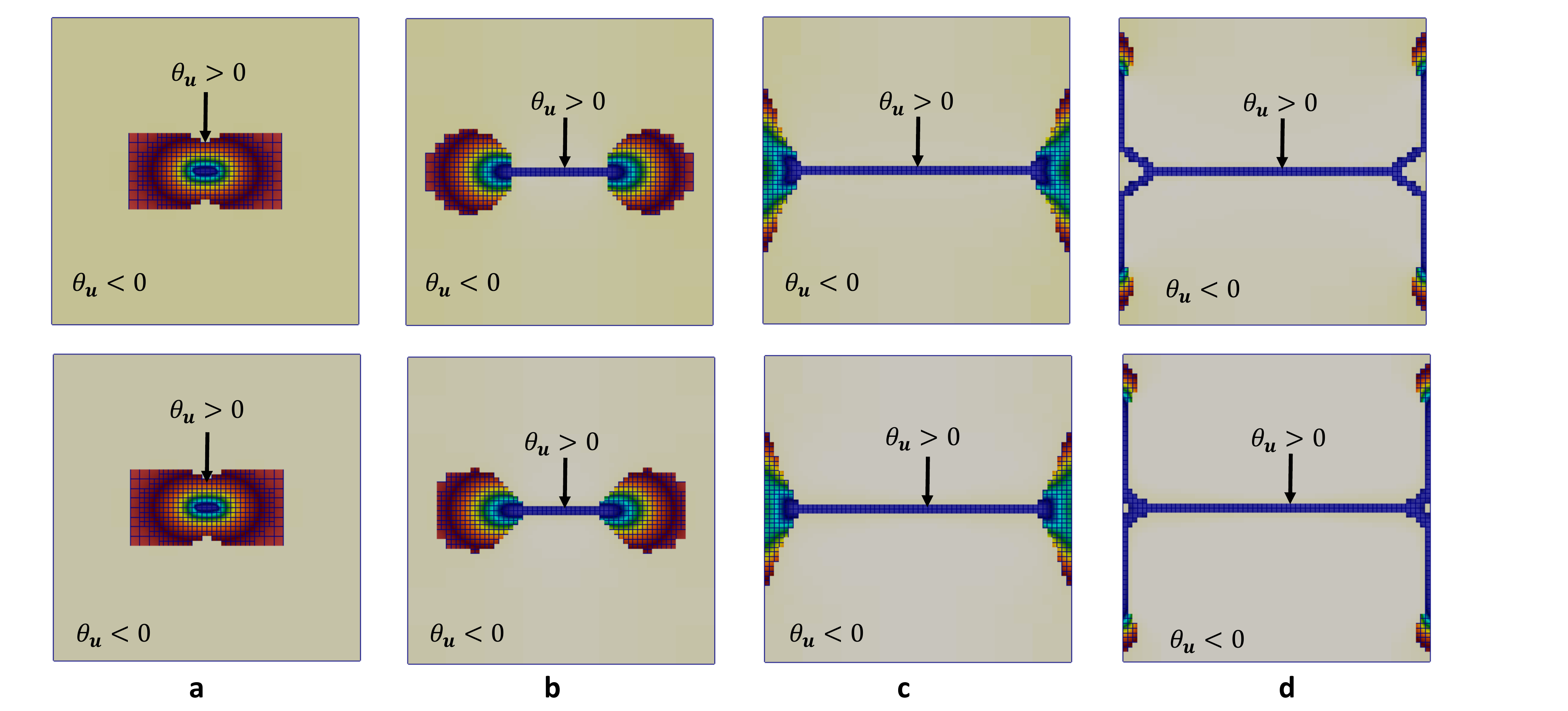}}   
	\caption{Case e. Qualitative representation of the 
tensile and compression counterpart of the domain at (a) 240 min (b) 252 min (c) 259 min (d)267 min.The first row represents no splitting and the second 
          row shows split of the strain energy density function. }
	\label{Figure15}
\end{figure}

A qualitative representation of the tensile region, i.e. $\theta_{\bm u}>0$, 
and compression counterpart, i.e. $\theta_{\bm u}<0$, are displayed in Fig. \ref{Figure15}.

Based on the results, shown in Fig. \ref{Figure15} (and also Fig. \ref{Figure12}), an 
alternative scalar variable for adaptive mesh refinement could be 
$\theta_{\bm u}:=H{^+}(\nabla. \bm u)$ rather than using the phase-field variable.
If $\theta_{\bm u}=1$, 
the material tends to tensile stress and hence crack propagation might 
occur and if $\theta_{\bm u}=0$ the material is under compression stress 
and we are away from the fracture zone. The main advantages of using  $\theta_{\bm u}$ 
rather than $\varphi\le \operatorname{\texttt{TOL}_\varphi}$ are
two-fold. First, the above suggestion is
based on the (physical) displacement field (and not on phase-field) and 
secondly, there is no requirement about choosing the correct 
threshold value for $\operatorname{\texttt{TOL}_\varphi}$.

\subsubsection{Performance of nonlinear and linear solvers}
\begin{figure}[!ht]
	\centering
	{\includegraphics[clip,trim=0cm 0cm 0cm 0cm, width=15cm] {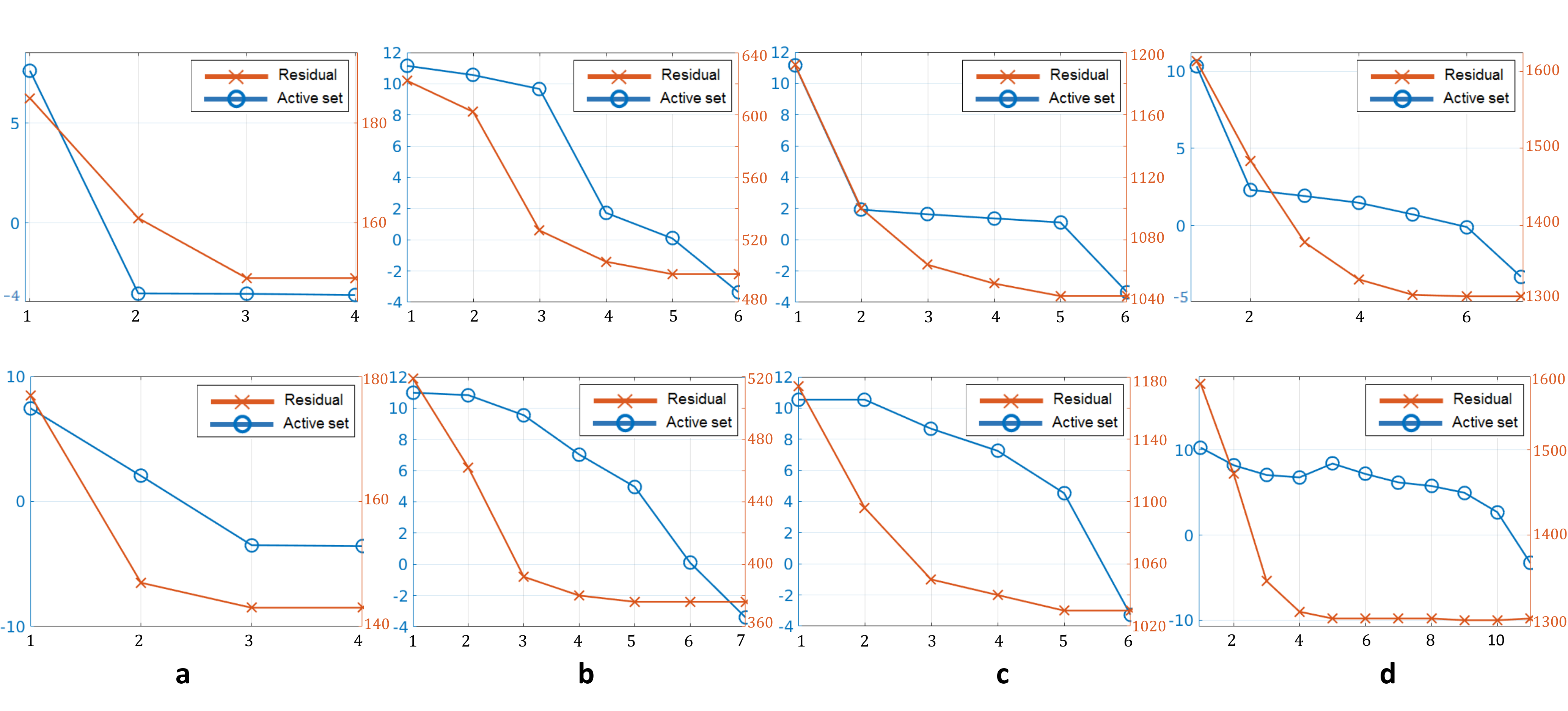}}   
	\caption{Case e. Performance of the combined Newton
          solver. Convergence behavior for the nonlinear relative residual,
          i.e. red line, and elements in the active set, i.e. blue line, at
          fixed time steps (a) 240 min (b) 252 min (c) 259 min (d) 267 min.
The first row represents no splitting and the second 
          row shows split of the strain energy density function. 
}
	\label{Figure151}
\end{figure}

As in Case d, we study the solver performances. The evolution of the 
primal-dual active set and the Newton residuals are shown in Fig. \ref{Figure151}. 

The average number of GMRES iterations of linear iterations per Newton cycle for 
different times steps is shown in Table \ref{GMRES_itr_caseE}. 
Specifically, the preconditioner (Section \ref{sec_lin_sol}) works again extremely well.

\begin{table}[!ht]\small
	\caption{Case e. Average number of GMRES iterations per Newton cycle for different times steps.}
	\label{GMRES_itr_caseE}
{
	\begin{center}
	\begin{tabular}{c|c|c}
		\hline 
		Time step (min) & no split & vol./dev. split \\
		\hline 
		240 & 2 & 2 \\ 
		252 & 10 & 14\\
		259 & 4 & 17\\
		267 & 13 & 24\\ 
		\hline
	\end{tabular}
	\end{center}
}
\end{table}

\newpage
\subsection{Case f. 3D setting with pressure, temperature and decline constant}
\begin{figure}[!ht]
	\centering
	{\includegraphics[clip,trim=0cm 2cm 0cm 2cm, width=16cm] {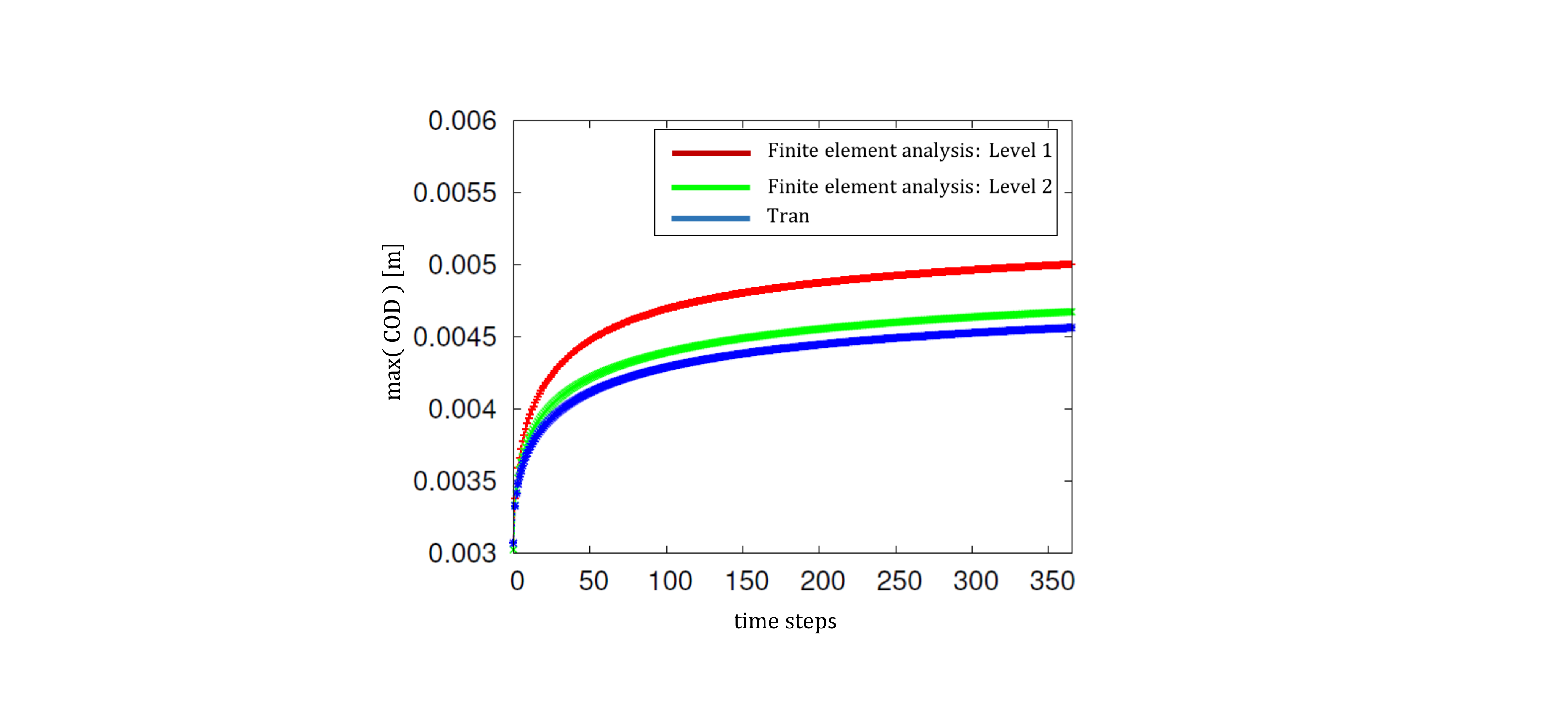}}   
	\caption{Case f. Maximum width evolution over time in the 
		non-isothermal case. Injection of cold water $\Theta < \Theta_0$ leads 
		to an increasing aperture over large time scales; here $365$ days.}
	\label{ex_1c}
\end{figure}

In this test, we extend Case c to a three-dimensional setting. The 
geometry and the material parameters are described in the Sections
\ref{sec_geometry_param} and \ref{sec_case_c}, respectively.
The maximum width evolution is displayed in Figure \ref{ex_1c}
and shows a good fit with Tran's et al. \cite{TrSeNg13} manufactured 3D
solution for computing the maximal COD.

\subsection{Case g. 3D setting for the Pressurized crack propagation in isothermal and non-isothermal settings}
\begin{figure}[!ht]
\centering
{\includegraphics[clip,trim=0cm 2cm 0cm 2cm, width=13cm]{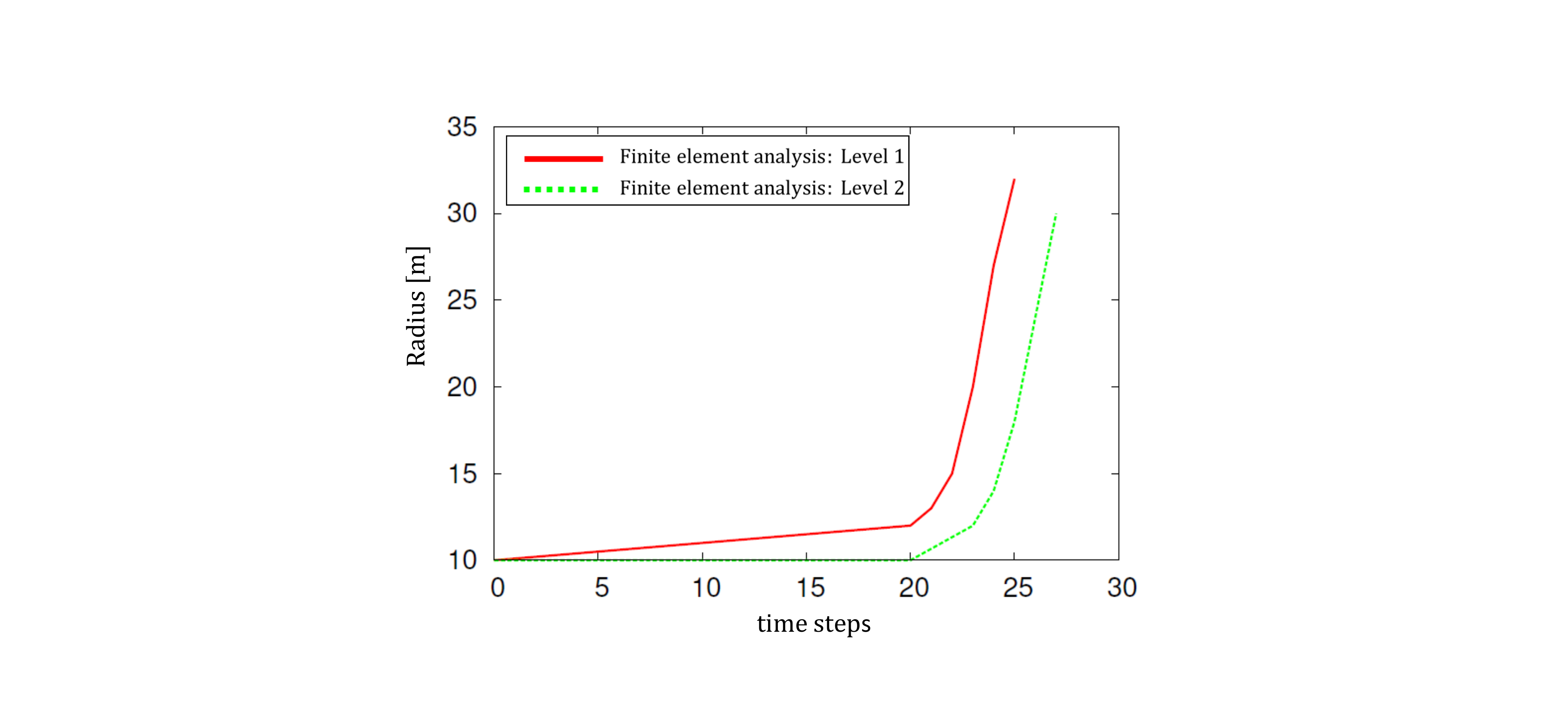}}
\caption{Case g. Radius evolution versus time
on two spatial mesh refinement levels. 
We observe that the crack on the finer mesh starts a bit later 
to grow but on both refinement levels the growth velocity (i.e., the slope)
is the same.}
\label{ex_3c}
\end{figure}

\begin{figure}[!ht]
	\centering
        {\includegraphics[clip,trim=0cm 0cm 0cm 0cm, width=14cm] {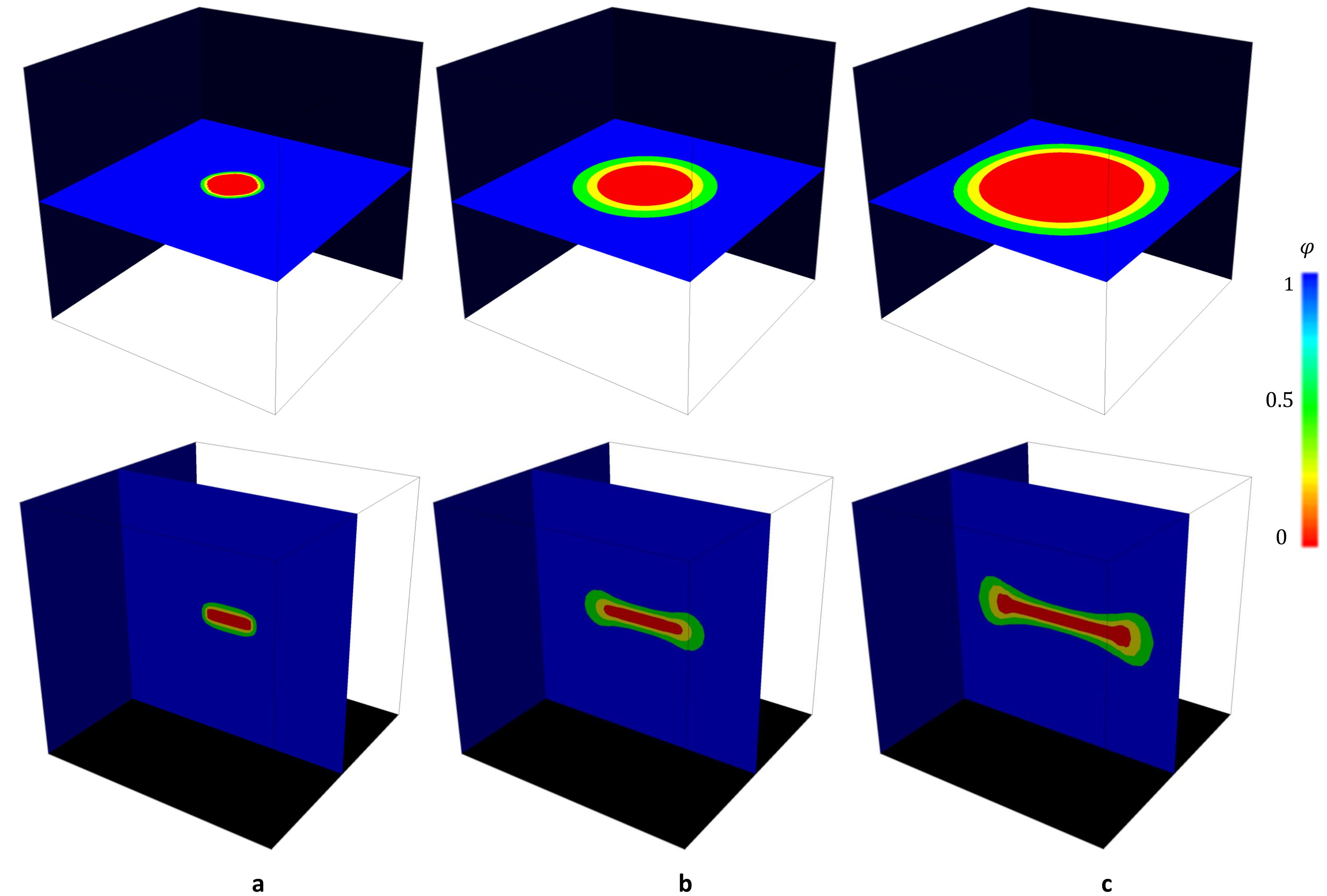}}
	\caption{Case g. Penny-shaped crack phase-field evolution at three
          times (a) $t=0$ min, (b) $t=25$ min	and (c) $t=27$ min.}
	\label{ex_2a}
\end{figure}

In this final example, the 2D Case e is extended to three dimensions.
The geometry and material parameters are listed 
in the Sections \ref{sec_geometry_param} and \ref{sec_case_e}.
As in the $2D$ test case, the time step size is $\delta t= 60s$.
We compute $27$ loading (time) steps. The fracture radii are displayed 
in Figure \ref{ex_3c}. Graphical illustrations of the phase-field variable,
i.e., the crack path, are displayed in Figure \ref{ex_2a}.

\section{Conclusions}
\label{sec_conc}
In this work, we developed a phase-field fracture 
model for pressurized and non-isothermal cracks in thermo-poroelasticity. 
We concentrated on the so-called mechanics step in which pressure 
and temperature are given quantities. Using the thermo-poroelastic 
stress tensor and interface laws, we derived an energy functional 
and the corresponding Euler-Lagrange equations. We then 
proceeded with a careful analysis from the mechanical point of view. 
As computational framework, we used a programming code 
that allows for adaptive mesh refinement and is fully parallelized.
Therein, the nonlinear problem is formulated within a monolithic framework and
is solved with a semi-smooth Newton method accounting as well 
for the crack irreversibility constraint.
Based on these computational capacities, we ran in total 
seven test scenarios. Therein, we compared to analytical (manufactured) 
solutions taken from the literature with an excellent agreement of 
our numerical solutions.
Moreover, we carried out mesh refinement studies and gave further 
insight into the performances of the linear and nonlinear 
solvers. In these studies, we also analyzed the effects on crack growth
and solver behavior when strain-energy splitting is used or not.
The outcome of our numerical simulations let us come 
to the conclusion that we have developed a robust and efficient 
computational framework for pressurized and non-isothermal 
configurations when the pressure and temperatures are given. In future studies, one idea is to 
treat the temperature as unknown to be determined by the governing partial
differential equation.
Here, we could proceed analogously to fluid-filled 
phase-field fracture models proposed in \cite{MiWheWi14}. 
Another task will be more computations of practical problems with multiple 
interacting fractures.

\section*{Acknowledgments}
This work is supported by the German Research Foundation, Priority Program 1748 (DFG SPP 1748) named
\textit{Reliable Simulation Techniques in Solid Mechanics. Development of
Non-standard Discretization Methods, Mechanical and Mathematical Analysis}
in the sub-project \textit{Structure Preserving Adaptive Enriched Galerkin Methods for Pressure-Driven 3D Fracture Phase-Field Models} (WI 4367/2-1).



\end{document}